\documentclass[sn-mathphys,Numbered]{sn-jnl}% Math and Physical Sciences Reference Style

\usepackage{amssymb,amsfonts,amsthm}%
\usepackage{lmodern}
\usepackage[title]{appendix}%
\usepackage{xcolor}%
\usepackage{manyfoot}%
\usepackage{booktabs}%
\usepackage{amsmath}
\usepackage{multirow}
\usepackage{tabularx}

\newcommand{\newsub}[1]{#1}
\newcommand{\martina}[1]{#1}

\begin{document}

\title[Article Title]{Neural networks for the approximation of Euler's elastica}

\author[1]{\fnm{Elena} \sur{Celledoni}}\email{elena.celledoni@ntnu.no}

\author[1]{\fnm{Ergys} \sur{\c{C}okaj}}\email{ergys.cokaj@ntnu.no}

% \author[1]{\fnm{Andrea} \sur{Leone}}\email{andrea.leone@ntnu.no}
\author[1]{\fnm{Andrea} \sur{Leone}}\email{andrea.leone@ntnu.no}

\author*[2]{\fnm{Sigrid} \sur{Leyendecker}}\email{sigrid.leyendecker@fau.de}

\author[1]{\fnm{Davide} \sur{Murari}}\email{davide.murari@ntnu.no}

\author[1]{\fnm{Brynjulf} \sur{Owren}}\email{brynjulf.owren@ntnu.no}

\author[2]{\fnm{Rodrigo T.} \sur{Sato Martín de Almagro}}\email{rodrigo.t.sato@fau.de}

\author[2]{\fnm{Martina} \sur{Stavole}}\email{martina.stavole@fau.de}

\affil*[1]{\orgdiv{Department of Mathematical Sciences}, \orgname{Norwegian University of Science and Technology}, \orgaddress{\street{Alfred Getz' vei 1}, \city{Trondheim}, \postcode{7034}, \country{Norway}}}

\affil[2]{\orgdiv{Institute of Applied Dynamics}, \orgname{Friedrich-Alexander-Universität Erlangen-Nürnberg}, \orgaddress{\street{Immerwahrstrasse 1}, \city{Erlangen}, \postcode{91058}, \country{Germany}}}

\abstract{Euler's elastica is a classical model of flexible slender structures relevant in many industrial applications. Static equilibrium equations can be derived via a variational principle. The accurate approximation of solutions to this problem can be challenging due to nonlinearity and constraints. We here present two neural network-based approaches for simulating Euler's elastica. Starting from a data set of solutions of the discretised static equilibria, we train the neural networks to produce solutions for unseen boundary conditions. We present a \textit{discrete} approach learning discrete solutions from the discrete data. We then consider a \textit{continuous} approach using the same training data set but learning continuous solutions to the problem. We present numerical evidence that the proposed neural networks can effectively approximate configurations of the planar Euler's elastica for a range of different boundary conditions.
}


\maketitle
\section{Introduction}\label{secintro}

Modelling of mechanical systems is relevant in various branches of engineering. Typically, it leads to the formulation of variational problems and differential equations, whose solutions are approximated with numerical techniques. The efficient solution of linear and nonlinear systems resulting from the discretisation of mechanical problems has been a persistent challenge of applied mathematics. While classical solvers are characterised by a well-established and mature body of literature \cite{saad2003iterative,nocedal1999numerical,marsden2001, hairer1993, hairer1996, brenner2008mathematical,quarteroni2006numerical}, the past decade has witnessed a surge in the use of novel machine learning-assisted techniques \cite{cuomo2022scientific,Brunton2023MachineLF, raissi2019physics, samaniego2020energy, yu2018deep, gu2023deep,kadupitiya2022solving,liu2020hierarchical,mattheakis2022hamiltonian,lu2021learning,lu2021deepxde,chevalier2022accelerating,li2022delisa,schiassi2021extreme,de2022physics,fabiani2023parsimonious,mortari2019high}. These approaches aim at enhancing solution methods by leveraging the wealth of available data and known physical principles. The use of deep learning techniques to improve the performance of traditional numerical algorithms in terms of efficiency, accuracy, and computational scalability \cite{Brunton2023MachineLF}, is becoming increasingly popular also in computational mechanics \newsub{(see, e.g. \cite{vuquoc2023})}. \newsub{Examples include a wide range of problems that require the approximation of functions, as well as efficient reduced order modeling \cite{brunton_kutz_2022} or more specific numerical tasks such as optimising the quadrature rule for computing the finite element stiffness matrix \cite{yagawa2022computational} or the investigation of data-driven numerical frameworks for the bifurcation analysis of partial differential equations \cite{fabiani2021numerical, galaris2022numerical}.} 
% \andout{Examples comprise virtually any problem where approximation of functions is required, but also efficient reduced order modelling e.g. in fluid mechanics, the deep Ritz method, or more specific numerical tasks such as optimisation of the quadrature rule for the computation of the finite element stiffness matrix, acceleration of simulations on coarser meshes by learning appropriate collocation points, and replacing expensive numerical computations with data-driven predictions \cite{Brunton2023MachineLF, kollmannsberger2021deep, yagawa2022computational, vuquoc2023,yu2018deep}.}
This recent literature is evidence that neural networks can be used successfully as surrogate models for the solution operators of various differential equations.

In the context of ordinary and partial differential equations, two main trends can be identified. The first one aims at providing a machine learning-based approximation to the discrete solutions of differential problems on a specific space-time grid, for example, by solving linear or nonlinear systems efficiently and accelerating convergence of iterative schemes \cite{li2022delisa, chevalier2022accelerating, liu2020hierarchical, kadupitiya2022solving, gu2023deep}. The second one provides instead solutions to the differential problem as continuous (and differentiable) functions of the temporal and spatial variables. Depending on the context, conditions on such approximate solutions are provided by the differential problem itself, the initial values and boundary conditions, and the available data. The idea of providing approximate solutions as functions defined on the space-time domain and parametrised as neural networks was proposed in the nineties \cite{lagaris1998artificial} and was recently revived in the framework of Physics-Informed Neural Networks in \cite{raissi2019physics}. Since then, such an approach has attracted much interest and developed in many directions \cite{cuomo2022scientific, schiassi2021extreme, kollmannsberger2021deep}.

In this work, we use neural networks to approximate the configurations of highly flexible slender structures modelled as beams. Such models are of great interest in industrial applications like cable car ropes, diverse types of wires or endoscopes \cite{ntarladima2023model,stavole2022homogenization,manfredo2023data,saadat2023mixed}. Notwithstanding their ingenious and simple mathematical formulation, slender structure models can accurately reproduce complex mechanical behaviour and for this reason their numerical discretisation is often challenging. Furthermore, the use of 3-dimensional models requires high computational time. Due to the fact that slender deformable structures have one dimension (length) being orders of magnitude larger than their other dimensions (cross-section), it is possible to reduce the complexity of the problem from a $3$-dimensional elastic continuum to a $1$-dimensional beam. A beam is modelled as a centerline curve, $\mathbf{q}: [0, L] \to \mathbb{R}^n, s \mapsto \mathbf{q}(s)$, with $n = 2$ or $n = 3$, along which a rigid cross-section $\Sigma(s)$ is attached. The main model assumption is that the diameter of $\Sigma(s)$ is small compared with the undeformed length $L$. The complexity of the model depends on factors such as the dimension of the problem, the translational and rotational degrees of freedom (DOF) at each node of the beam, the analysis, i.e., static or dynamic. Exploring the numerous beam models documented in the literature, we choose to approach the challenge of approximating beam deformations using a simple yet widely employed model, i.e., the $2$-dimensional \textit{Euler's elastica} \cite{Euler}. The cross-section $\Sigma(s)$ 
is assumed to have unchanged geometrical and material properties, and be orthogonal to the centerline $\mathbf{q}(s)$. The latter is an inextensible curve and solution of a bending energy minimisation problem \cite{love1863, matsutani2010,singer2008} for given boundary conditions.

Although the 2-dimensional Euler's elastica is relatively simple compared to more comprehensive models, it can robustly represent interesting real-world phenomena. For instance, the elastica model appropriately captures the high bending deformations of flexible endoscopes, complex medical devices, during surgeries \cite{stavole2022homogenization}. The approximation of the elastica through neural networks can help predict the deformed configuration of the beam for endoscopy simulations, particularly when the beam encounters constraints in confined spaces.

When approximating static equilibria of Euler's elastica via neural networks, a key issue is to ensure the inextensibility of the curve (having unit norm tangents) as well as the boundary conditions. Two main approaches can be found in the literature \cite{kollmannsberger2021deep, schiassi2021extreme, rohrhofer2022on}. One is the weak imposition of constraints and boundary conditions adding appropriate extra terms to the loss function. The other is a strong imposition strategy consisting in shaping the network architectures to satisfy the constraints by construction. We show examples of both the approaches in Sections~\ref{nndiscrete} and \ref{sec:nncont}.

The paper is organised as follows. In Section \ref{secmodel}, we present the mathematical model of the planar Euler's elastica, including its continuous and discrete equilibrium equations. We describe the approach used to generate the data sets for the numerical experiments. In Section \ref{secdeeplearning}, we introduce some basic theory and notation for neural networks that we shall use in the succeeding sections. Starting from general theory, we specialise in the task of approximating configurations of Euler's elastica. In Section \ref{nndiscrete}, we introduce the \textit{discrete} approach, which aims to approximate precomputed numerical discretisations of Euler's elastica. \newsub{This represents the natural approach to approximate the discrete solution trajectories with a parametric method.} We discuss some drawbacks associated with this approach and then propose an alternative approximation strategy in Section \ref{sec:nncont}, \newsub{that leverages the fact that we are approximating a continuous curve on a spatial grid}. The \textit{continuous} approach consists in computing an arc length parametrisation of the beam configuration.
We provide insights into two additional networks and analyse how the test accuracy changes with varying constraints, such as boundary conditions or tangent vector norms. Data and codes for the numerical experiments are available in the GitHub repository associated to the paper\footnote{\href{https://github.com/ergyscokaj/LearningEulersElastica}{https://github.com/ergyscokaj/LearningEulersElastica}}.\\ 

\textbf{\textit{Main contributions:}} This paper presents advancements in the approximation of beam \newsub{static} configurations using neural networks. These advancements include: (i) A \newsub{detailed} experimental analysis of approximating numerical discretisations of Euler's elastica configurations through what we call \textit{discrete network}, (ii) Identification and discussion of the limitations associated with this discrete approach, and (iii) Introduction of a new parametrisation strategy called \textit{continuous network} to address some of these drawbacks.
\begin{table*}[ht]
\centering
\small
{
\begin{tabularx}{1.\textwidth}{ 
   >{\centering\arraybackslash}p{0.2\textwidth}
   >{\arraybackslash}p{0.75\textwidth}}
\toprule
    \multicolumn{2}{ c }{Nomenclature}\\
    \midrule
    $\mathcal{L}$ & \text{continuous Lagrangian function}\\
    $\mathcal{S}$ & \text{continuous action functional}\\
    $\mathcal{L}_d$ & \text{discrete Lagrangian function}\\
    $\mathcal{S}_d$ & \text{discrete action functional}\\
    $\mathbf{q}$ & configuration of the beam \\
    $\mathbf{q}'$ & first spatial derivative of $\mathbf{q}$ \\
    $\theta$ & tangential angle\\
    $s$ & arc length parameter \\
    $\kappa$ & curvature \\
    %$\mathcal{P}$ & \text{nonlinear differential operator related to the problem} \\ 
    $L$ & length of the undeformed beam \\
    $EI$ & bending stiffness, with $E$ the elastic modulus and $I$ the second moment of area \\
    $\hat{\mathbf{q}}$ & numerical approximation of $\mathbf{q}$ \\
    $N+1$ & number of discretisation nodes, with $N$ the number of intervals\\
    $h$ & \text{space step (length of each interval)}\\
    %\text{NN} & \text{neural network}\\
    $q_{\boldsymbol{\rho}}^{\textrm{d}}$ & \text{discrete neural network}\\  $q_{\boldsymbol{\rho}}^{\textrm{c}}$ & \text{continuous neural network approximating the solution curve $\boldsymbol{q}(s)$}\\ $\theta_{\boldsymbol{\rho}}^{\textrm{c}}$ & \text{continuous neural network approximating the angular function $\theta(s)$}\\
    $\boldsymbol{\rho}$ & \text{parameters of the neural network}\\
    $\ell$ & number of layers in the neural network \\
    $\sigma$ & \text{activation function} \\
    $M$ & number of training data \\
    $B$ & size of one training batch \\
    \text{MSE} & \text{mean squared error}\\
    \text{MLP} & \text{multi layer perceptron}\\
    \text{MULT} & \text{multiplicative neural network}\\
    $\mathcal{D}$ & \text{differential operator}\\
    $\mathcal{I}$ & \text{quadrature operator}\\
    \bottomrule
\end{tabularx}
\caption{List of abbreviations and notations.}
\label{tab:nomenclature}}
\end{table*}
\section{Euler's elastica model}\label{secmodel}
We consider an inextensible beam model in which the cross-section $\Sigma(s)$ is assumed to be constant along the arc length $s$ and perpendicular to the centerline $\mathbf{q}(s)$, which means that no shear deformation can occur. Thus, the deformation of the centerline is a pure bending problem, precisely Euler's elastica curve. In the following, we assume $\mathbf{q} \in C^2([0,L], \mathbb{R}^2)$, i.e., the curve is planar and twice continuously differentiable with length $L$. If $s$ denotes the arc length parameter, then $\|\mathbf{q}'(s)\| = 1$, where $'=\frac{d}{ds}$, for all $s \in [0,L]$. The elastica problem consists in minimising the following Euler-Bernoulli energy functional
{
\begin{equation*}
    \int_0^L \kappa(s)^2 ds,
\end{equation*}
}
where $\kappa(s)$ denotes the curvature of $\mathbf{q}(s)$, \cite{matsutani2010}. Given the arc length parametrisation, then $\kappa(s) = \|\mathbf{q}''(s)\|$.\\

We can reformulate this problem as a constrained Lagrangian problem as follows. Consider the second-order Lagrangian $\mathcal{L}: T^{(2)}Q \rightarrow \mathbb{R}$, where $T^{(2)}Q$ denotes the second-order tangent bundle \cite{colombo2016} of the configuration manifold $Q$, which in this case is $\mathbb{R}^2$:
{
\begin{equation}\label{eq:L_unconstrained_elastica}
   \mathcal{L}\left(\mathbf{q}, \mathbf{q}', \mathbf{q}''\right) = \frac{1}{2} EI \left\|\mathbf{q}''\right\|^2\,.
\end{equation}
}
Here, abusing the notation, ${}^\prime$ denotes a spatial derivative, but we do not initially assume arc length parametrisation. The parameter $EI$ is the bending stiffness, which governs the response of the elastica under bending. This mechanical parameter consists of a material and a geometric properties, where $E$ is the Young's modulus and $I$ is the second moment of area of the cross-section $\Sigma$. For simplicity, these parameters are assumed to be constant along the length of the beam.\\

In order to recover the solutions of the elastica, the Lagrangian in Equation~\eqref{eq:L_unconstrained_elastica} must be supplemented with the constraint equation
{
\begin{equation} \label{eq:constraint}
    \Phi(\mathbf{q},\mathbf{q}')=\|\mathbf{q}'\|^2 - 1 = 0.
\end{equation}
}
This imposes arc length parametrisation of the curve $\mathbf{q}(s)$ and leads to the augmented Lagrangian $\widetilde{\mathcal{L}}: T^{(2)}Q \times \mathbb{R} \rightarrow \mathbb{R}$
{
\begin{equation}\label{eq:L_static_elastica}
    \widetilde{\mathcal{L}} \left(\mathbf{q}, \mathbf{q}', \mathbf{q}'',\Lambda\right) = \mathcal{L}\left(\mathbf{q}, \mathbf{q}', \mathbf{q}''\right) + \Lambda \Phi(\mathbf{q},\mathbf{q}'),
\end{equation}
}
where $\Lambda(s)$ is a Lagrange multiplier, see \cite{singer2008}. The Lagrangian function coincides with the total elastic energy over solutions of the corresponding Euler-Lagrange equations. The internal bending moment is directly related to the curvature $\kappa(s)$.\\

The continuous action functional $\mathcal{S}$ is defined as:
{
\begin{equation}
    \label{eq:action_integral}
    \mathcal{S}[\mathbf{q}] = \int_{0}^{L} \widetilde{\mathcal{L}} \left(\mathbf{q},\mathbf{q}',\mathbf{q}'',\Lambda \right) \,ds.
\end{equation}
}
Applying Hamilton's principle of stationary action, $\delta \mathcal{S} = 0$, yields the Euler-Lagrange equations 
{
\begin{equation}\label{eq:second_order_ELeq}
    \begin{aligned}
				\frac{d^2}{d s^2} \left( \frac{\partial \mathcal{L}}{\partial \mathbf{q}''} \right) - \frac{d}{d s} \left( \frac{\partial \mathcal{L}}{\partial \mathbf{q}'} \right) + \frac{\partial \mathcal{L}}{\partial \mathbf{q}} &= \frac{d}{d s} \left( \frac{\partial \Phi}{\partial \mathbf{q}'} \Lambda \right) - \frac{\partial \Phi}{\partial \mathbf{q}} \Lambda,\\
				\|\mathbf{q}'\|^2 - 1 &= 0,
	\end{aligned}
\end{equation}
}
which need to be satisfied together with the boundary conditions on positions and tangents, i.e., $(\mathbf{q}(0),\mathbf{q}'(0)) = (\mathbf{q}_0,\mathbf{q}'_0)$ and $(\mathbf{q}(L),\mathbf{q}'(L)) = (\mathbf{q}_N, \mathbf{q}'_N)$.

\subsection{Space discretisation of the elastica}\label{eqsdiscrete}

The continuous augmented Lagrangian $\widetilde{\mathcal{L}}$ in Equation~\eqref{eq:L_static_elastica} and the action integral $\mathcal{S}$ in Equation~\eqref{eq:action_integral} are discretised over the beam length $L$ \newsub{using constant step size $h=L/N$, with $N+1$ the number of the resulting equidistant nodes $0=s_0<s_1<\ldots<s_{N-1}<s_N=L$}. In second-order systems, the discrete Lagrangian is a function $\widetilde{\mathcal{L}}_d: TQ \times TQ \times \mathbb{R} \times \mathbb{R} \rightarrow \mathbb{R}$. In this study, we refer to a discretisation of the Lagrangian function proposed in \cite{ferraro2021} based on the trapezoidal rule:

{
\[
\begin{split}
&\widetilde{\mathcal{L}}_d \left( \mathbf{q}_k, \mathbf{q}'_k, \mathbf{q}_{k+1}, \mathbf{q}'_{k+1}, \Lambda_k, \Lambda_{k+1} \right)\\
&=
    \frac{h}{2} \left[ \widetilde{\mathcal{L}} (\mathbf{q}_k,\mathbf{q}'_k,\left(\mathbf{q}''_k\right)^{-},\Lambda_k) + \widetilde{\mathcal{L}} (\mathbf{q}_{k+1},\mathbf{q}'_{k+1},\left(\mathbf{q}''_{k+1}\right)^{+}, \Lambda_{k+1}) \right],
\end{split}
\]
}
where $\mathbf{q}_k$, $\mathbf{q}'_k$, and $\Lambda_k$ are approximations of $\mathbf{q}(s_k)$, $\mathbf{q}'(s_k)$, and $\Lambda(s_k)$, and the curvature on the interval $[s_k, s_{k+1}]$ is approximated in terms of lower order derivatives as follows 
{
    \begin{align*}
    \mathbf{q}''(s_k)\approx (\mathbf{q}''_k)^{-} =& \, \frac{\left( -2 \mathbf{q}'_{k+1} - 4 \mathbf{q}'_k \right) h + 6  (\mathbf{q}_{k+1} - \mathbf{q}_k)}{h^2}\,, \\
    \mathbf{q}''(s_{k+1})\approx(\mathbf{q}''_{k+1})^{+} =& \, \frac{\left( 4 \mathbf{q}'_{k+1} + 2 \mathbf{q}'_k \right) h - 6 (\mathbf{q}_{k+1} - \mathbf{q}_k)}{h^2}\,.
\end{align*}
}
This amounts to a piece-wise linear and discontinuous approximation of the curvature on $[0,L]$.

The action integral in Equation~\eqref{eq:action_integral} along the exact solution $\mathbf{q}$ with boundary conditions $\left( \mathbf{q}_0,\mathbf{q}'_0 \right)$ and $\left( \mathbf{q}_N,\mathbf{q}'_N \right)$ is approximated by 
{
\begin{equation} \label{eq:discrete_action}
    \mathcal{S}_d = \sum_{k=0}^{N-1} \widetilde{\mathcal{L}}_d \left( \mathbf{q}_k, \mathbf{q}'_k, \mathbf{q}_{k+1}, \mathbf{q}'_{k+1}, \Lambda_k, \Lambda_{k+1} \right).
\end{equation}
}
The discrete variational principle $\delta \mathcal{S}_d=0$ leads to the following discrete Euler-Lagrange equations:
{
\begin{equation}\label{eq:DEL_beam}
    \begin{aligned}
			D_3 \widetilde{\mathcal{L}}_d \left( \mathbf{q}_{k-1},\mathbf{q}'_{k-1},\mathbf{q}_k,\mathbf{q}'_k,\Lambda_{k-1}, \Lambda_{k} \right) + D_1 \widetilde{\mathcal{L}}_d \left( \mathbf{q}_k,\mathbf{q}'_k,\mathbf{q}_{k+1},\mathbf{q}'_{k+1},\Lambda_{k},\Lambda_{k+1} \right)  &= 0, \\
			D_4 \widetilde{\mathcal{L}}_d \left( \mathbf{q}_{k-1},\mathbf{q}'_{k-1},\mathbf{q}_k,\mathbf{q}'_k,\Lambda_{k-1},\Lambda_{k} \right) + D_2 \widetilde{\mathcal{L}}_d \left( \mathbf{q}_k,\mathbf{q}'_k,\mathbf{q}_{k+1},\mathbf{q}'_{k+1},\Lambda_{k},\Lambda_{k+1} \right) &= 0, \\
			D_6 \widetilde{\mathcal{L}}_d \left( \mathbf{q}_{k-1},\mathbf{q}'_{k-1},\mathbf{q}_k,\mathbf{q}'_k,\Lambda_{k-1},\Lambda_{k} \right)  + D_5 \widetilde{\mathcal{L}}_d \left( \mathbf{q}_k,\mathbf{q}'_k,\mathbf{q}_{k+1},\mathbf{q}'_{k+1},\Lambda_{k},\Lambda_{k+1} \right) &= 0,
		\end{aligned}
\end{equation}
}
for $k=1,\dots,N-1$, which approximate the equilibrium equations of the beam in Equations~\eqref{eq:second_order_ELeq} and can be solved together with the boundary conditions. \newsub{Here, $D_i$ for $i = 1, \dots, 6$ denotes the differentiation with respect to the $i$-th argument.}

\subsection{Data generation} \label{subsec:bifurcation}
The elastica was one of the first examples displaying elastic instability and bifurcation phenomena \cite{timoshenko1961,bigoni2012}. Elastic instability implies that small perturbations of the boundary conditions might lead to large changes in the beam configuration, which results in unstable equilibria. Under certain boundary conditions, bifurcation can appear leading to a multiplicity of solutions \cite{matsutani2010}. In particular, this means that the numerical problem may display history-dependence and converge to solutions that do not minimise the bending energy. In order to generate a physically meaningful data set, avoiding unstable and non-unique solutions is essential. Thus, in addition to the minimisation of the discrete action $S_d$ in Equation~\eqref{eq:discrete_action}, we ensure the fulfilment of the discrete Euler-Lagrange equations \eqref{eq:DEL_beam}, which can be seen as necessary conditions for the stationarity of the discrete action. We exclude from the data set numerical solutions computed with boundary conditions where minimisation of Equation~\eqref{eq:discrete_action} and accurate solution of Equations~\eqref{eq:DEL_beam} can not be simultaneously achieved. 

In particular, we consider a curve of length $L=3.3$ and bending stiffness $EI=10$, divided into $N = 50$ intervals. We fix the endpoints $\mathbf{q}_0 = (0,0)$, $\mathbf{q}_N=(3,0)$. The units of measurement are deliberately omitted as they have no impact on the results of this work. We impose boundary conditions on the tangents in the following two variants:
\begin{enumerate}
    \item the angle of the tangents with respect to the $x$-axis at the boundary, $\theta_0$ and $\theta_N$, is prescribed in the range $[0, 2 \pi]$, in a specular symmetric fashion, i.e., $\theta_N = \pi -\theta_0$. Hereafter, we refer to this case as \textit{both-ends},
   \item the angle of the left tangent is left fixed as $\theta_0=0$ and the angle of the right tangent, $\theta_N$, varies in the range of $[0,2\pi]$. We refer to this case as \textit{right-end}.
\end{enumerate}

Based on these parameters and boundary values, \newsub{and using cubic splines as initial guess,} we generate a data set of $2000$ trajectories ($1000$ trajectories for each case) by minimising the particular action in Equation~\eqref{eq:discrete_action}, with the \texttt{trust-constr} solver of the \texttt{optimize.minimize} procedure provided in \texttt{SciPy} \cite{2020SciPy-NMeth}. We check the resulting solutions by using them as initial guesses for the \texttt{optimize.root} method of \texttt{SciPy}, solving the discrete Euler-Lagrange equations~\eqref{eq:DEL_beam}.

\newsub{The learning problem we consider relies on numerically generated solution curves. This choice allows us to work with data points that are quantifiably close to the analytical solution of Euler's elastica. Consequently, showing that the neural networks we propose can accurately approximate these curves translates into their ability to approximate the analytical solution accurately. The motivation of this strategy is not to improve on the numerical solver we use, but to use its accuracy to train a model able to extrapolate to unseen boundary conditions and generate their solution curves more efficiently than the numerical method itself. The chosen supervised learning setting is independent of the fact that we use numerical solutions as data. Indeed, if one had another reliable approximation of the analytical solution, for example, based on realistic measurements, those could also be used or combined with numerically generated trajectories. Using numerical solutions as data is not an inherent limitation of the proposed procedure but a choice we make to quotient out the issue of not having reliable input data.} \martina{Furthermore, we mainly focus on the development of neural networks able to approximate such input data with high accuracy.}

\section{Approximation with neural networks}\label{secdeeplearning}
\newsub{We start by providing a concise overview of neural networks, which also serves to define the notation used in Sections \ref{nndiscrete} and \ref{sec:nncont}. We refer to \cite{higham2019deep, kollmannsberger2021deep,Goodfellow-et-al-2016} and references therein for a more extensive introduction.} A neural network is a parametric function $f_{\boldsymbol{\rho}}:\mathcal{I}\rightarrow\mathcal{O}$ with parameters $\boldsymbol{\rho} \in \Psi$ given as a composition of multiple transformations,
{
\begin{equation} \label{eq:fnncomposition}
f_{\boldsymbol{\rho}} := f_{\ell} \circ\dots\circ f_j \circ \dots \circ f_1,
\end{equation}
}
where each $f_j$ represents the $j$-th layer of the network, with $j=1,\dots,\ell$, and $\ell$ is the number of layers. For example, multi-layer perceptrons (MLPs) have each layer $f_j$ defined as
{
\begin{equation} \label{eq:mlpnn}
    f_j^{\mathrm{MLP}}(\mathbf{x}) = \sigma(\mathbf{A}_j\mathbf{x}+\mathbf{b}_j)\in\mathbb{R}^{n_j},
\end{equation}
}
where \newsub{$n_j$ is the dimension of the output of the $j$-th layer,} $\mathbf{x}\in \mathbb{R}^{n_{j-1}}$, and $\mathbf{A}_j \in \mathbb{R}^{n_j \times n_{j-1}} $, $\mathbf{b}_j\in \mathbb{R}^{n_j}$ are the parameters of the $j$-th layer, i.e., $\boldsymbol{\rho} = \{\mathbf{A}_j, \mathbf{b}_j\}_{j=1}^{\ell}$. The activation function $\sigma$ is a continuous nonlinear scalar function, which acts component-wise on vectors. The architecture of the neural network is prescribed by the layers $f_j$ in Equation~\eqref{eq:fnncomposition} and determines the space of functions $\mathcal{F}=\{f_{\boldsymbol{\rho}}:\mathcal{I}\rightarrow\mathcal{O},\,\,\boldsymbol{\rho}\in\Psi\}$ that can be represented. The weights $\boldsymbol{\rho}$ are chosen such that $f_{\boldsymbol{\rho}}$ approximates accurately enough a map of interest $f:\mathcal{I}\rightarrow\mathcal{O}$. Usually, this choice follows from minimising a purposely designed loss function $\rm{Loss}(\boldsymbol{\rho})$. 

In supervised learning, we are given a data set $\Omega=\{\mathbf{x}^i, \mathbf{y}^i\}_{i=1}^{M}$ consisting of $M$ pairs $(\mathbf{x}^i,\mathbf{y}^i=f\left(\mathbf{x}^i\right))$. The loss function measures the distance between the network predictions $f_{\boldsymbol{\rho}}\left(\mathbf{x}^i\right)$ and the desired outputs $\mathbf{y}^i$ in some appropriate norm $\|\cdot \|$,
{
\begin{equation*} \label{eq:mseloss}
\textrm{Loss}(\boldsymbol{\rho}) = \frac{1}{M}\sum_{i=1}^{M} \left\|f_{\boldsymbol{\rho}}\left(\mathbf{x}^i\right) - \mathbf{y}^i \right\|^2.
\end{equation*}
}
The training of the network is the process of minimising $\rm{Loss}(\boldsymbol{\rho})$ with respect to $\boldsymbol{\rho}$ and it is usually done with gradient descent (GD):
{
\begin{equation*}\label{eq:gd}
\boldsymbol{\rho}{^{(k)}} \mapsto \boldsymbol{\rho}{^{(k)}}-\eta \nabla \textrm{Loss}\left(\boldsymbol{\rho}{^{(k)}}\right)=:{\boldsymbol{\rho}^{(k+1)}}.
\end{equation*}
}
The scalar value $\eta$ is known as the learning rate. The iteration process is often implemented using subsets of data $\mathcal{B}\subset\Omega$ of cardinality $B=|\mathcal{B}|$ (batches). In this paper we use an accelerated version of GD known as Adam \cite{Adam_KingBa15}.

\newsub{During training, we evaluate the model's prediction accuracy using inputs in a validation set. This helps to prevent overfitting on the training data and may serve as a stopping criterion if the training loss diminishes but the validation error rises. Once the training is complete, we assess the model's accuracy in predicting the correct output for new inputs included in a test set composed of boundary conditions outside the training and validation sets. In the following, we measure the accuracy on the training, validation, and test data using the mean squared error of the difference between the predicted trajectories and the true ones.}

We now turn to the task of approximating the static equilibria of the planar elastica introduced in Section~\ref{secmodel}, i.e., approximating a family of curves $\{\mathbf{q}^i:[0,L] \mapsto\mathbb{R}^2 \}$ determined by boundary conditions,
{
\begin{equation} \label{bdata}
\{\mathbf{q}^i(0)=\mathbf{q}^i_0, \; \mathbf{q}^i(L)=\mathbf{q}^i_N, \; (\mathbf{q}^i)'(0)=(\mathbf{q}^i_0)', \; (\mathbf{q}^i)'(L)=(\mathbf{q}^i_N)'\},
\end{equation}
}
where $\left(\mathbf{q}^i_0,\mathbf{q}^i_N,(\mathbf{q}^i_0)',(\mathbf{q}^i_N)'\right)\in\mathbb{R}^8$. To tackle this problem, we require %either 
a set of evaluations $\{\mathbf{q}^i_k, (\mathbf{q}^i_k)'\}$ on the nodes $s_k \in [0,L]$ of a discretisation. More precisely, in our setting, the data set includes numerical approximations $\hat{\mathbf{q}}$ of the solution $\mathbf{q}(s)$ and its spatial derivative $\mathbf{q}'(s)$ at the $N-1$ discrete locations $s_k = \frac{kh}{L}$ in the interval $[0,L]$, for $M$ pairs of boundary conditions, as described in Section~\ref{subsec:bifurcation}.

\section{The discrete network}\label{nndiscrete}
The discretisation of Euler's elastica presented in Section \ref{eqsdiscrete} provides discrete solutions on a set of nodes along the curve. These solutions can sometimes be hard to obtain since a \newsub{global} optimisation problem needs to be solved, and the number of nodes can be large. This motivates using neural networks to learn the approximate solution on the internal nodes for a given set of boundary conditions. The data set $\Omega$ consists of $M$ precomputed discrete solutions
{
$$
\Omega = \left\{(\mathbf{x}^i,\mathbf{y}^i)\right\}_{i=1}^M,
$$
}
where 
{
\begin{equation*}
    \mathbf{x}^i=\left(\mathbf{q}^i_0,(\mathbf{q}^i_0)',\mathbf{q}^i_N,(\mathbf{q}^i_N)'\right)\in\mathbb{R}^8
\end{equation*}
}
are the input boundary conditions and
{
\begin{equation*}
    \mathbf{y}^i=(\hat{\mathbf{q}}_1^i,(\hat{\mathbf{q}}^{i}_1)',\ldots,\hat{\mathbf{q}}_{N-1}^i, (\hat{\mathbf{q}}^{i}_{N-1})')\in\mathbb{R}^{4(N-1)}
\end{equation*}
}
are the computed solutions at the internal nodes that serve as output data for the network's training. 

For any symmetric positive definite matrix $W$, we define the weighted norm $\|\mathbf{x}\|_{W}^2=\mathbf{x}^{\top}W\mathbf{x}$.
The weighted MSE loss
{
\begin{equation} \label{lossdiscrete}
    \mathrm{Loss}(\boldsymbol{\rho}) = \frac{1}{4M(N-1)}\sum_{i=1}^M\left\|q_{\boldsymbol{\rho}}^{\textrm{d}}\left(\mathbf{x}^i\right)- \mathbf{y}^i\right\|_{W}^2
\end{equation}
}
will be used to learn the input-to-output map $q_{\boldsymbol{\rho}}^{\textrm{d}}: \mathbb{R}^8 \to \mathbb{R}^{4(N-1)}$, where the superscript $\mathrm{d}$ stands for discrete. One should be aware that there is a numerical error in $\mathbf{y}^{i}$ compared to the exact solution and the size of this error will pose a limit to the accuracy of the neural network approximation.

\subsection{Numerical experiments}
\begin{figure}[ht]
\centering
\includegraphics[width=0.48\textwidth]{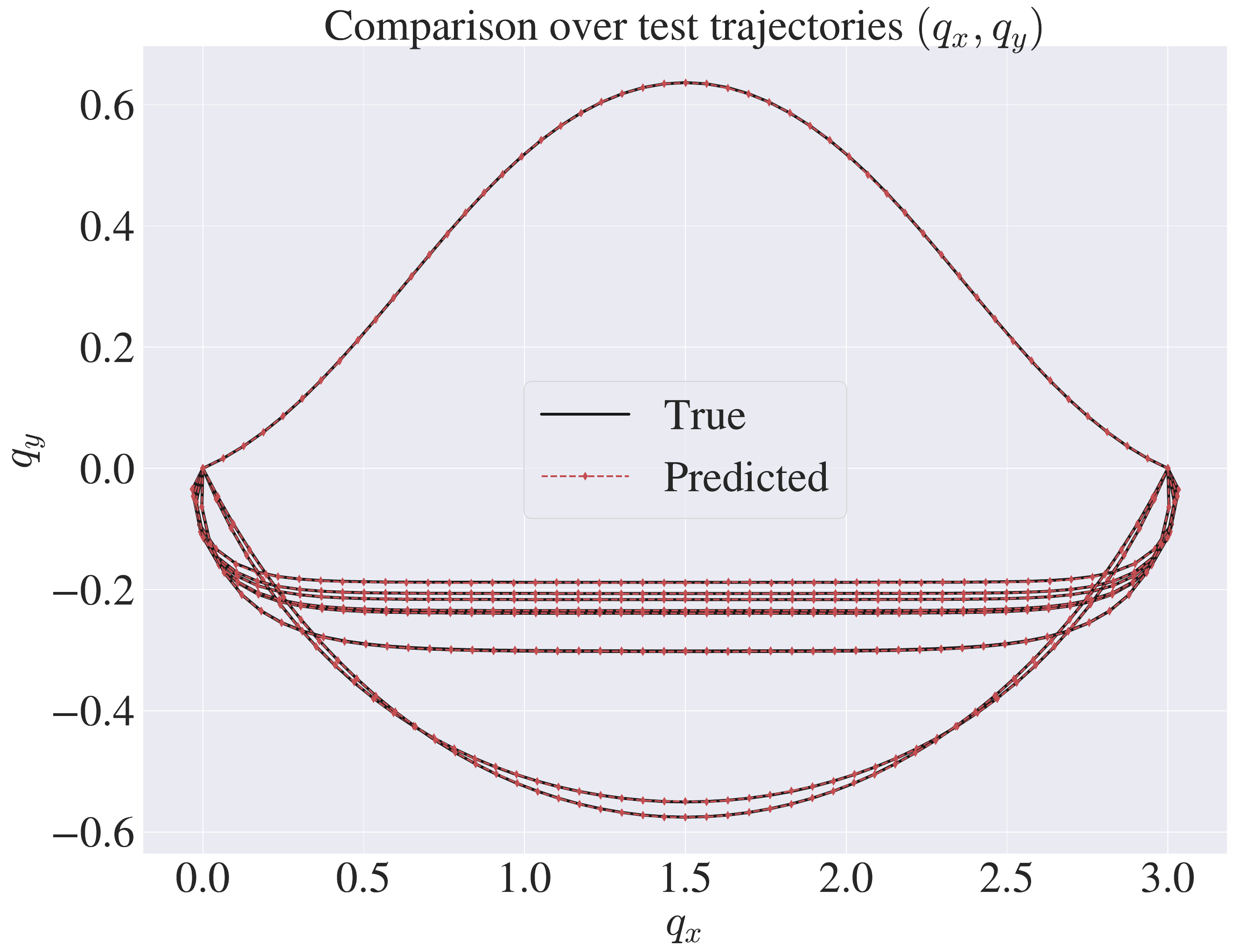}
\includegraphics[width=0.49\textwidth]{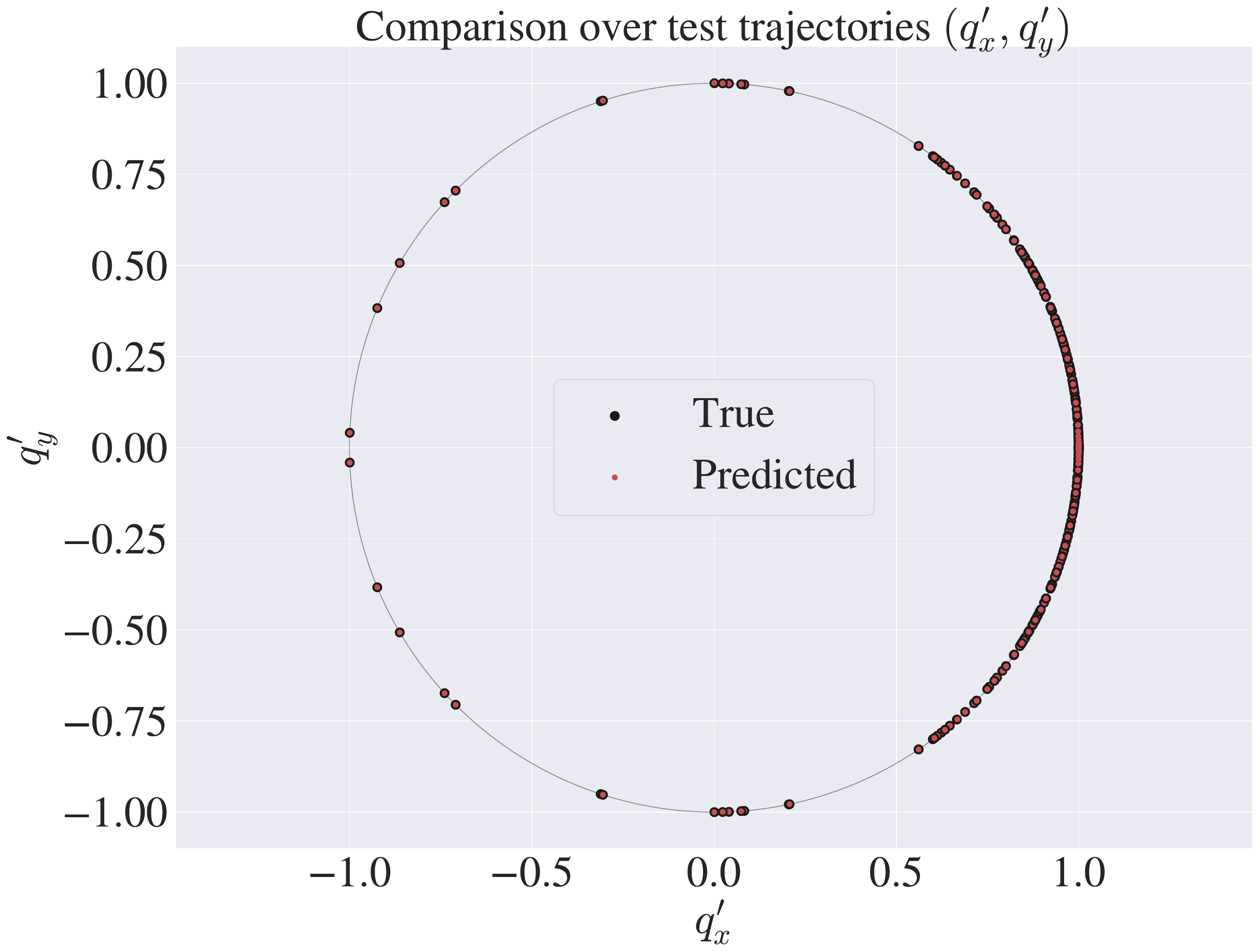}
\includegraphics[width = 0.49\textwidth]{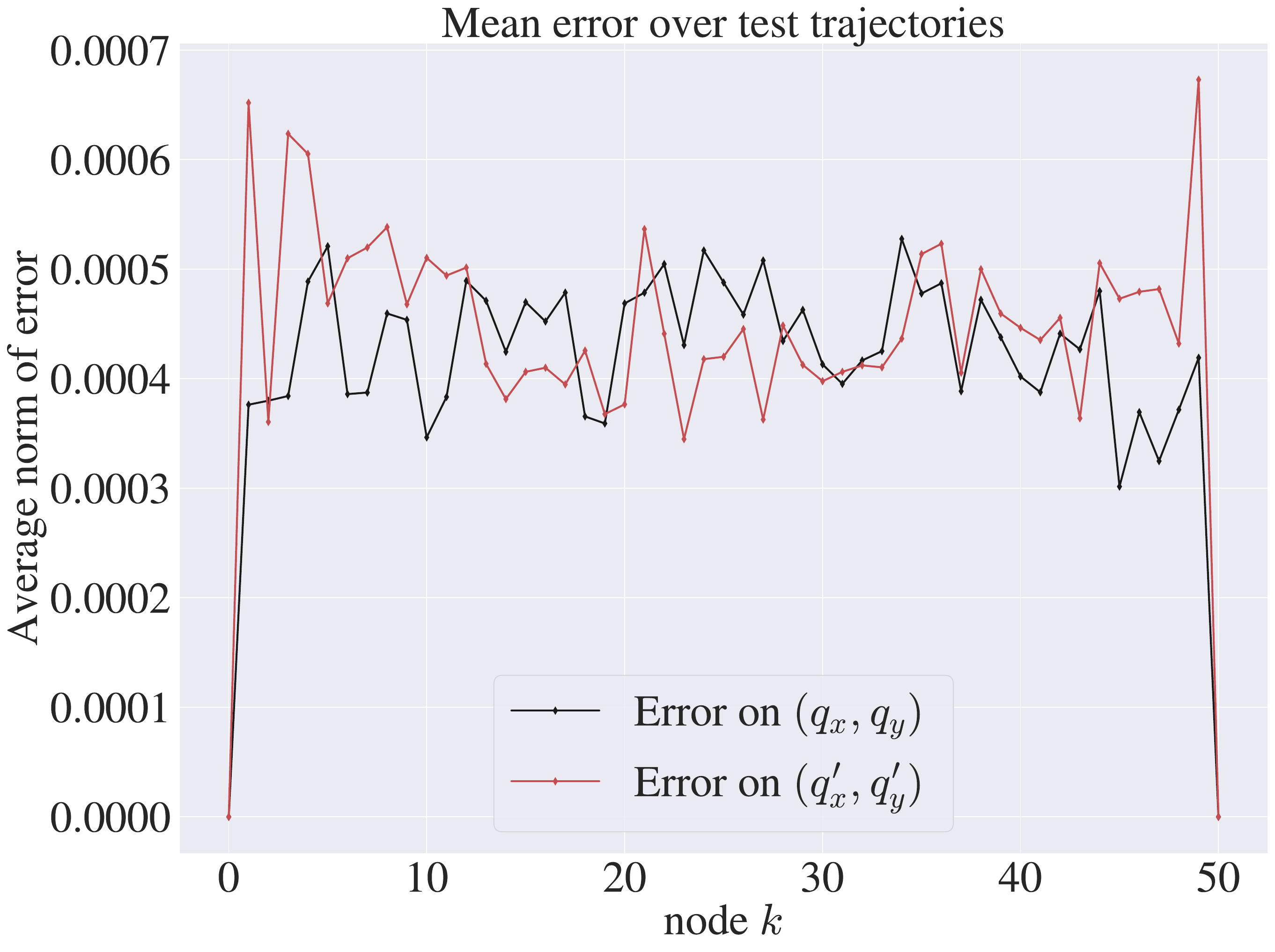}
\caption{Comparison over test trajectories for $\mathbf{q}$ and $\mathbf{q}^{\prime}$ for the discrete network $q_{\boldsymbol{\rho}}^{\textrm{d}}$ tested on the \textit{both-ends} data set with $80\% \-- 10\%\-- 10\%$ splitting into training, validation, and test sets. 
% \andout{These results are obtained with the hyperparameters from Table \ref{tab:comp_table_hyperparameters_discrete_net}, that yield a training error equal to $1.14 \cdot 10^{-7}$, a validation error equal to $2.151\cdot 10^{-7}$, and a test error equal to $4.009 \cdot 10^{-7}$.} 
The mean squared error on the test set equals $4.009 \cdot 10^{-7}$. For presentation purposes, only 10 randomly selected trajectories are considered in the first two plots.}
\label{fig:q_qp_discrete_net}
\end{figure}

This section provides experimental support to the proposed learning framework \newsub{using the machine learning library \texttt{PyTorch} \cite{paszke2019pytorch}.} \newsub{The experiments of this section are run on a CPU machine.} We perform a series of experiments \newsub{varying some hyperparameters in the training procedure. We fix the batch size $B$ to 32
% , the total number of epochs to $300$, 
and use the Adam optimiser \cite{Adam_KingBa15} for the training with learning rate $10^{-3}$ and weight decay set to $0$.} In \eqref{lossdiscrete} we use the weight matrix
{
$$W=I+\gamma G^\top G,$$
}
where $G=S^4-I$ with $S$ the forward shift operator on vectors of $\mathbb{R}^{4(N-1)}$. This choice of $G$ allows us to compute differences between corresponding entries of the input associated with neighbouring nodes. \newsub{We determine the number of epochs for training both the discrete and continuous networks based on experimental evidence. We fix a high enough number which allows us to achieve qualitatively accurate predictions and ensure that both training and validation losses start to plateau a few epochs before the set maximum.} \newsub{We consider a multi-layer perceptron with the hyperbolic tangent as an activation function, and we vary the number of layers and the number of hidden nodes in each layer. We also test different values of the parameter $\gamma$ in the weight matrix $W$.} We rely on the software framework \texttt{Optuna} \cite{akiba2019optuna}, \newsub{which employs Bayesian optimisation methods to automate and efficiently conduct the search for the combination that yields the best result.} We collect in Table~\ref{tab:hyperparams_dicrete_net} the hyperparameters with the corresponding ranges and in Table~\ref{tab:comp_table_hyperparameters_discrete_net} the selected values. The resulting training error on the \textit{both-end} data set is $1.14 \cdot 10^{-7}$, the validation error is $2.151\cdot 10^{-7}$, and the test error is $4.009 \cdot 10^{-7}$. Figure~\ref{fig:q_qp_discrete_net} compares test trajectories for $\mathbf{q}$ and $\mathbf{q}^{\prime}$. We remark that, as already clear from the low value of the training and test errors, the network can accurately replicate the behaviour of the training and test data. Furthermore, we have zero errors at the end nodes since the network is trained only on the internal nodes and the boundary values are appended to the predicted solution in a post-processing phase. On the other hand, since this discrete approach does not relate the components as evaluations of a smooth curve, there is no regular behaviour in the error.

\begin{figure}[ht]
\centering
\includegraphics[width=0.48\textwidth]{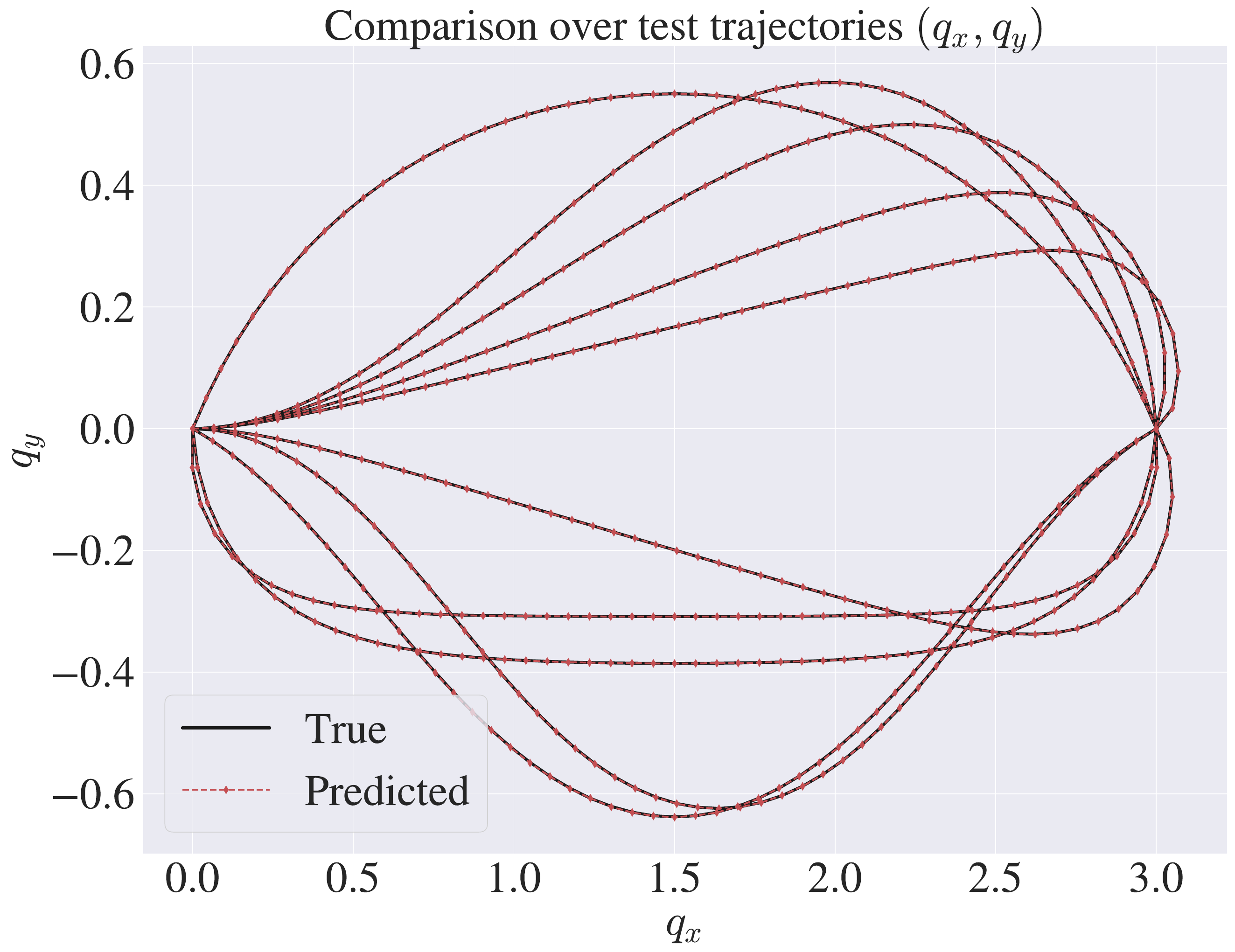}
\includegraphics[width=0.49\textwidth]{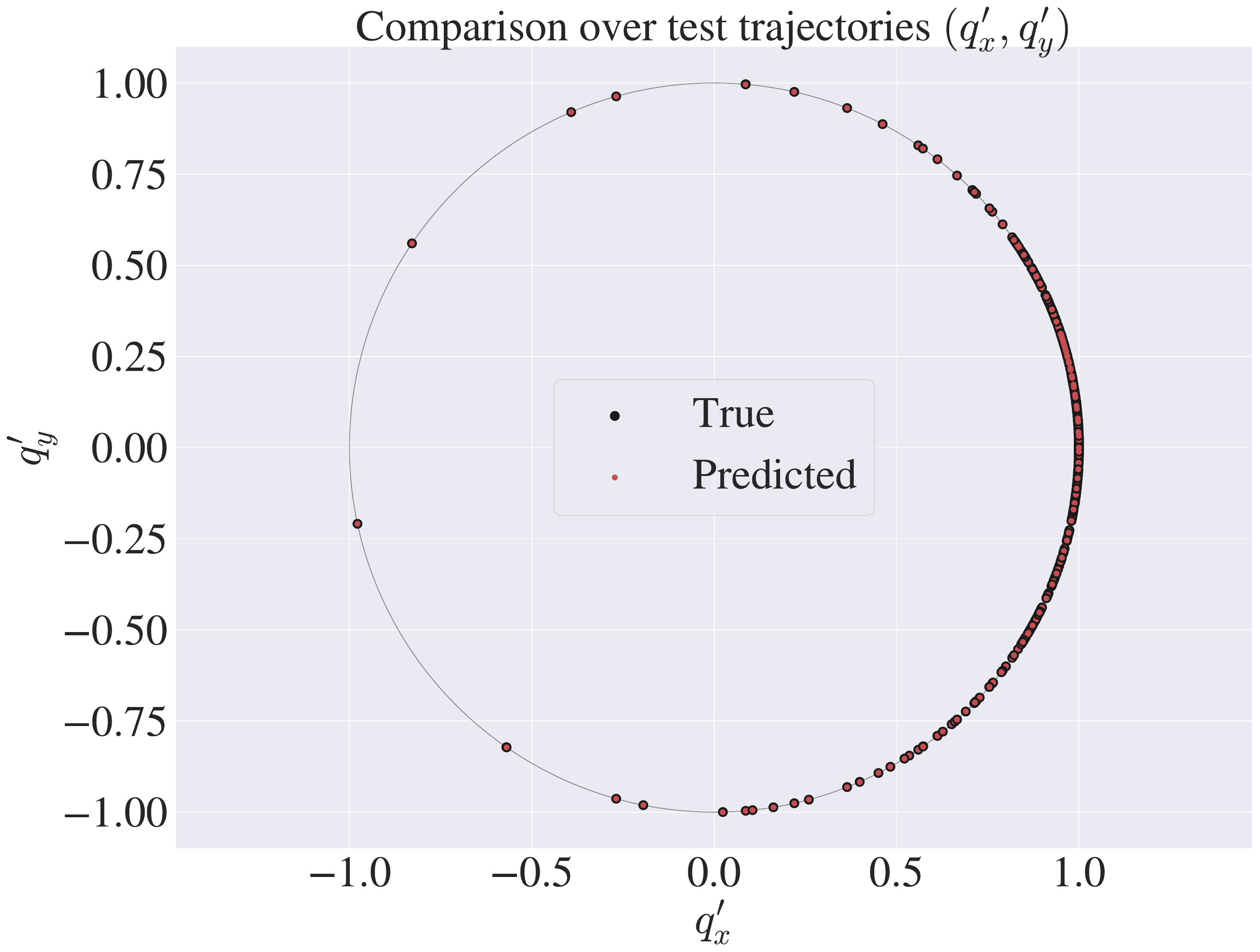}
\vfill
\includegraphics[width=0.49\textwidth]{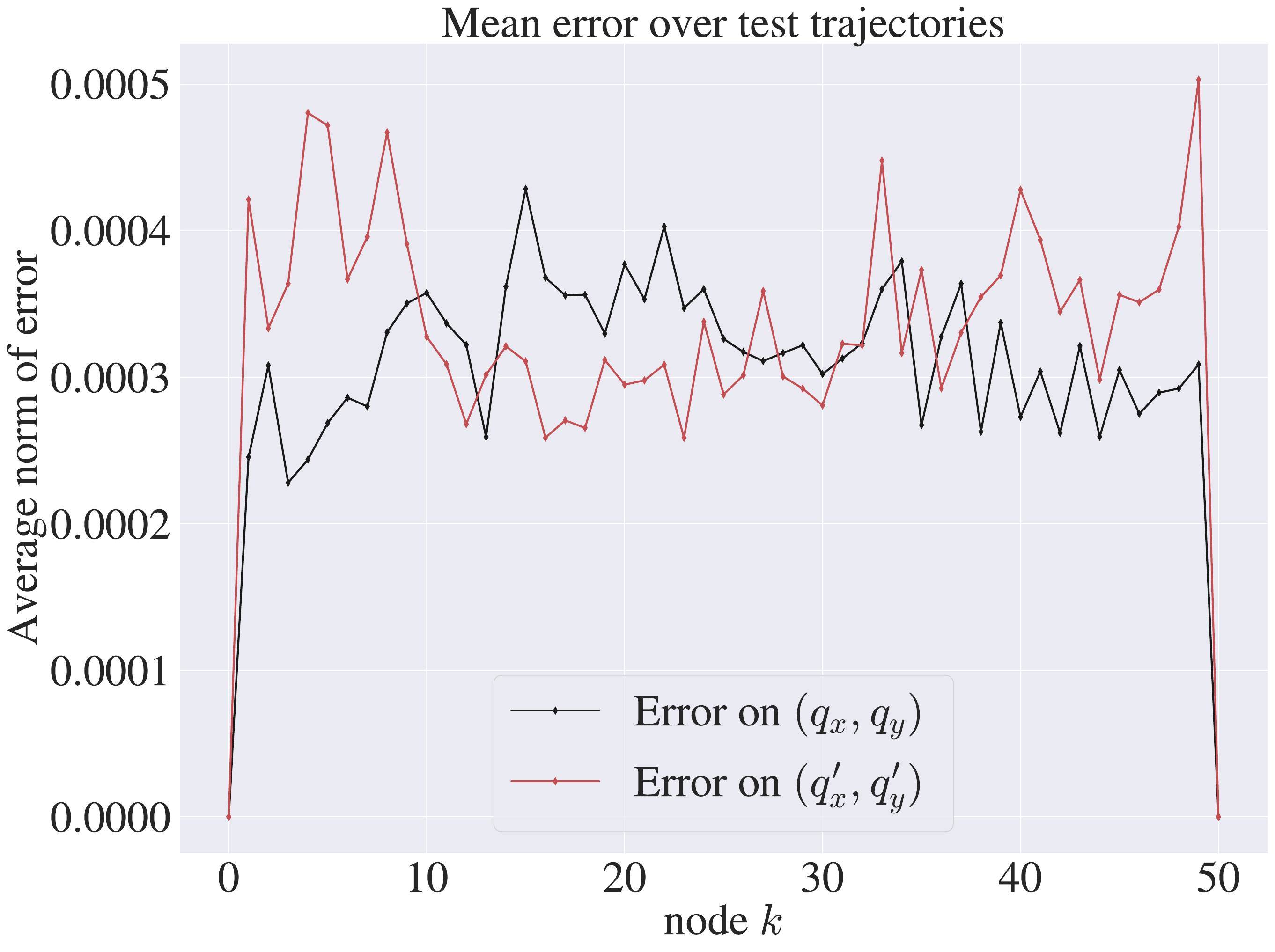}
\caption{Comparison over test trajectories for $\mathbf{q}$ and $\mathbf{q}^{\prime}$ for the discrete network $q_{\boldsymbol{\rho}}^{\textrm{d}}$ tested on the \textit{both-ends} + \textit{right-end} data set with $80\% \-- 10\%\-- 10\%$ splitting into training, validation, and test sets. 
% \ergys{\sout{These results are obtained with 3 layers, 616 hidden nodes, and $\gamma = 7.323 \cdot 10^{-3}$, that yield a training error equal to $9.893 \cdot 10^{-8}$, a validation error equal to $1.126\cdot 10^{-7}$, and a test error equal to  $ 7.854 \cdot 10^{-8}$.}} 
The mean squared error on the test set equals $ 7.854 \cdot 10^{-8}$.
For presentation purposes, only 10 randomly selected trajectories are considered in the first two plots.
}
\label{fig:error_discrete_net_both_right_end}
\end{figure}

As an additional evaluation of the deep learning framework's behaviour, we conduct experiments to assess how the learning process performs when the number of training data varies, i.e., with different splittings of the data set into training, validation, and test sets. We report the results in Table~\ref{tab:comp_table_accuracy_discrete_net} and summarise the corresponding hyperparameters in Table~\ref{tab:comp_table_hyperparameters_discrete_net} of the Appendix.
\begin{table*}[ht]
\small
\begin{tabularx}{1.\textwidth}{ 
   >{\centering\arraybackslash}p{0.28\textwidth} 
   >{\centering\arraybackslash}p{0.21\textwidth}   
   >{\centering\arraybackslash}p{0.21\textwidth}
   >{\centering\arraybackslash}p{0.21\textwidth}}
\toprule
$\begin{array}{c}
    \text{Data set splitting} \\
    \text{Training \-- validation \-- test}
    \end{array}$ & \text {Training accuracy} & \text {Validation accuracy} & \text{Test accuracy} \\
    \midrule \text {10\% \-- 10\% \-- 10\%} & $   2.331\cdot 10^{-5}$ & $   3.874\cdot 10^{-5}$ & $   8.545\cdot 10^{-4}$\\
    \text {20\% \-- 10\% \-- 10\%} & $  1.852 \cdot 10^{-6}$  & $   1.327\cdot 10^{-6}$ & $   1.361\cdot 10^{-4}$\\
    \text {40\% \-- 10\% \-- 10\%} & $   4.802\cdot 10^{-7}$ & $   4.793\cdot 10^{-7}$ & $  1.295\cdot 10^{-6}$ \\
    \text {80\% \-- 10\% \-- 10\%} & $  1.140\cdot 10^{-7}$ & $  2.151\cdot 10^{-7}$ & $   4.009\cdot 10^{-7}$\\
    \bottomrule
\end{tabularx}\caption{Behaviour of the discrete network $q_{\boldsymbol{\rho}}^{\textrm{d}}$ tested on the \textit{both-ends} data set with fewer training data points. The size of the training set varies, while that of the validation and the test sets is fixed.
The last row corresponds to the results in Figure \ref{fig:q_qp_discrete_net}.
}
\label{tab:comp_table_accuracy_discrete_net} 
\end{table*}
We also report results obtained by merging the \textit{both-end} and the \textit{right-end} trajectories, with $80\% \-- 10\%\-- 10\%$ splitting of the whole new data set into training, validation, and test sets. The results are shown in Figure \ref{fig:error_discrete_net_both_right_end} and are obtained with 3 layers, 616 hidden nodes, and $\gamma = 7.323 \cdot 10^{-3}$. The resulting training, validation, and test errors are, respectively, 
$9.893 \cdot 10^{-8}$, $1.126\cdot 10^{-7}$, and $ 7.854 \cdot 10^{-8}$.

\section{The continuous network} \label{sec:nncont}

The approach described in the previous section shows accurate results, given a large enough amount of beam discretisations with a fixed number of nodes $N+1$, equally distributed in $[0,L]$. It seems reasonable to expect the parametric model's approximation quality to improve when the number of discretisation nodes increases. However, in this approach, the dimension of the predicted vector grows with $N$, and hence minimising the loss function \eqref{lossdiscrete} becomes more difficult. In addition, the fact that the discrete network approach depends on the spatial discretisation of the training data restricts the output dimension to a specific number of nodes. Consequently, there would be two main options to assess the solution at different locations: training the network once more, or interpolating the previously obtained approximation. These limitations make such a discrete approach less appealing and suggest that having a neural network that is a smooth function of the arc length coordinate $s$ can be beneficial. This modelling assumption would also be compatible with different discretisations of the curve and would not suffer from the curse of dimensionality if more nodes were added.
In this setting, the discrete node $s_k$ at which an approximation of the solution is available, is included in the input data together with the boundary conditions. As a result, we work with the following data set   
{
$$
\Omega = \left\{\left(s_k,\,\mathbf{x}^i\right),\; {\mathbf{y}}_k^i \right\}_{k=0,\dots,N}^{i=1,\dots,M},
$$
where, as in the previous section, 
\[
\mathbf{x}^i = \left(\,\mathbf{q}^i_0,\,(\mathbf{q}^i_0)',\,\mathbf{q}^i_N,\,(\mathbf{q}^i_N)'\right)\in\mathbb{R}^8,
\]
and
\[
{\mathbf{y}}_k^i = \left(\hat{\mathbf{q}}_k^i, \, (\hat{\mathbf{q}}_k^i)' \right).
\]
}
Here $\hat{\mathbf{q}}_k^i$ is the numerical solution $\hat{\mathbf{q}}$ on the node $s_k$, satisfying the $i$-th boundary conditions in Equation~\eqref{bdata}. Let us introduce the neural network
\[
q_{\boldsymbol{\rho}}^{\mathrm{c}}:\mathbb{R}^8 \to \mathcal{C}^{\infty}\left([0,L],\mathbb{R}^2\right),
\]
and the differential operator
\[
\mathcal{D} : \mathcal{C}^{\infty}\left([0,L],\mathbb{R}^2\right)\to \mathcal{C}^{\infty}\left([0,L],\mathbb{R}^2\right),\,\,\mathcal{D}\left(q_{\boldsymbol{\rho}}^{\mathrm{c}}\left(\mathbf x^i\right)\right)(s_k) = \frac{d}{ds}\left(q_{\boldsymbol{\rho}}^{\rm c}(\mathbf{x}^i)\right)(s)\Big\vert_{s=s_k},
\]
%\partial_s q_{\boldsymbol{\rho}}^{\mathrm{c}}\left(\mathbf x^i\right)(s)\vert_{s=s_k}
so that we can define
\[
y_{\boldsymbol{\rho}}\left(\mathbf{x}^i\right)(s_k) := \left(q_{\boldsymbol{\rho}}^{\mathrm{c}}\left(\mathbf x^i\right)(s_k),\, \mathcal{D}\left(q_{\boldsymbol{\rho}}^{\mathrm{c}}\left(\mathbf x^i\right)\right)(s_k)\right).
\]
To train the network $q_{\boldsymbol{\rho}}^{\mathrm{c}}$, we define the loss function
{
\begin{equation}\label{eq:lossContinuos}
\begin{aligned}
    \textrm{Loss}(\boldsymbol{\rho})=& \frac{1}{4M(N+1)}\sum_{i=1}^M\sum_{k=0}^{N}\left( \left\|{ y_{\boldsymbol{\rho}}}\left(\mathbf x^i\right)(s_k) - {\boldsymbol{y}}_k^i\right\|^2_2 + \right. \\
    &\left.\gamma\left(\left\|\pi_{\mathcal{D}}\left(y_{\boldsymbol{\rho}}\left(\mathbf x^i\right)(s_k)\right)\right\|_2^2 - 1\right)^2 \right),
\end{aligned}
\end{equation}
}
where $\pi_{\mathcal{D}}:\mathbb{R}^8\to\mathbb{R}^4$ is the projection on the second component $\mathcal{D}(q_{\boldsymbol{\rho}}^{\rm c}(\mathbf{x}^i))(s_k)$, and $\gamma\geq 0$ weighs the violation of the normality constraint. The map $q_{\boldsymbol{\rho}}^{\textrm{c}}$ is now a neural network that associates each set of boundary conditions $\mathbf{x}^i$ with a smooth curve $q_{\boldsymbol{\rho}}^{\textrm{c}}\left(\mathbf{x}^i\right):[0,L]\to \mathbb{R}^2$ that can be evaluated at every point $s\in [0,L]$. We denote this network with the superscript $\mathrm{c}$ since this curve is, in particular, continuous. The outputs $q_{\boldsymbol{\rho}}^{\textrm{c}}\left(\mathbf{x}^i\right)(s)\in \mathbb{R}^2$ are approximations of the configuration of the beam at $s \in [0,L]$. 

We point out that, contrary to the discrete case, we learn approximations of $\mathbf{q}(s)$ also on the end nodes, i.e., at $s=0$ and $s=L$. This is because we do not impose the boundary conditions by construction. Even though there are multiple approaches to embed them into the network architecture, the one we try in our experiments made the optimisation problem too complex, thus we only impose the boundary conditions weakly in the loss function.

Another strategy is to compute the angles $\theta_{k}$ between the tangents $(\hat{\mathbf{q}}_k)'$ and the $x$-axis and to use them as training data. To this end,
we define the neural network 
\[
\theta_{\boldsymbol{\rho}}^{\textrm{c}} : \mathbb{R}^8 \to \mathcal{C}^{\infty}\left([0,L],\mathbb{R}\right)
\]
as $\theta_{\boldsymbol{\rho}}^{\textrm{c}} = \hat{\theta}_{\boldsymbol{\rho}}^{\textrm{c}} \circ \pi$, where 
\begin{equation}\label{eq:thetaHat}
\hat{\theta}_{\boldsymbol{\rho}}^{\textrm{c}} : \mathbb{R}^2 \to \mathcal{C}^{\infty}\left([0,L],\mathbb{R}\right)
\end{equation}
is a neural network, and the function $\pi:\mathbb{R}^8\to\mathbb{R}^2$ extracts the tangential angles from the boundary conditions, i.e., $\pi\left(\mathbf{x}^i\right) = \left(\theta_0^i,\theta_N^i\right)$. Such a network should approximate the angular function $\theta: [0,L] \ni s \to \mathbb{R}$, so that 
\begin{equation}\label{eq:tau_rec}
\tau_{\boldsymbol{\rho}}^{\mathrm{c}}\left(\mathbf{x}^i\right)(s):=\left(\cos\left(\theta_{\boldsymbol{\rho}}^{\textrm{c}}\left(\mathbf{x}^i\right)(s)\right),\sin\left(\theta_{\boldsymbol{\rho}}^{\textrm{c}}\left(\mathbf{x}^i\right)(s)\right)\right)\in\mathbb{R}^2
\end{equation}
gets close to the tangent vector $\mathbf{q}'(s)$. As a result, the constraint on the unit norm of the tangents is satisfied by construction, and the inextensibility of the elastica is guaranteed. The curve 
{
\begin{equation*}\label{eq:globalGauss}
\mathbf{q}(s) = \mathbf{q}_0 + \int_0^s \mathbf{q}'(\bar{s})\mathrm{d}\bar{s}
\end{equation*}
}
can then be approximated through the reconstruction formula 
{
\begin{equation}\label{eq:q_rec}
q_{\boldsymbol{\rho}}^{\textrm{c}}\left(\mathbf{x}^i\right)(s) = \mathbf{q}_0 + \mathcal{I}\left(\tau_{\rho}^{\mathrm{c}}\left(\mathbf{x}^i\right)\right)(s),
\end{equation}
}
where the operator $\mathcal{I} : \mathcal{C}^{\infty}\left([0,L],\mathbb{R}^2\right)\to \mathcal{C}^{\infty}\left([0,L],\mathbb{R}^2\right)$ is such that
\[
\mathcal{I}\left(\tau_{\rho}^{\mathrm{c}}\left(\mathbf{x}^i\right)\right)(s) \approx\int_0^{s}\tau_{\rho}^{\mathrm{c}}\left(\mathbf x^i\right)(\bar{s})\mathrm{d}\bar{s}.
\]
In the numerical experiments, $\mathcal{I}$ is based on the $3$-point Gaussian quadrature formula applied to a partition of the interval $[0,L]$, see \cite[Chapter 9]{quarteroni2006numerical}.
As done previously, we define the vector
\begin{equation}\label{eq:yTheta}
y_{\boldsymbol{\rho}}\left(\mathbf{x}^i\right)(s_k) := \left(q_{\boldsymbol{\rho}}^{\mathrm{c}}\left(\mathbf x^i\right)(s_k),\, \tau_{\boldsymbol{\rho}}^{\rm c}\left(\mathbf{x}^i\right)(s_k)\right),
\end{equation}
with components defined as in Equations \eqref{eq:tau_rec} and \eqref{eq:q_rec}. This allows us to train the network $\theta_{\boldsymbol{\rho}}^{\rm c}$ by minimising the same loss function as in Equation \eqref{eq:lossContinuos}, where this time $y_{\boldsymbol{\rho}}^{\rm c}$ is given by Equation \eqref{eq:yTheta}. Furthermore, since by construction this case satisfies $\left\|\pi_{\mathcal{D}}\left(y_{\boldsymbol{\rho}}^{\rm c }(\mathbf{x}^i)(s)\right)\right\|_2 = \left\|\tau_{\boldsymbol{\rho}}^{\mathrm{c}}\left(\mathbf{x}^i\right)(s)\right\|_2\equiv 1$, we set $\gamma=0$. We present numerical experiments for the two proposed continuous networks $q_{\boldsymbol{\rho}}^{\textrm{c}}$ and $\theta_{\boldsymbol{\rho}}^{\textrm{c}}$. In the latter case, by neural network architecture, we refer to $\hat{\theta}_{\boldsymbol{\rho}}^{\rm c}$ rather than $\theta_{\boldsymbol{\rho}}^{\rm c}$ in what follows. We analyse $q_{\boldsymbol{\rho}}^{\textrm{c}}$ more thoroughly in Section~\ref{subsec:num_exp_r_rho_c}, mirroring most of the discrete case experiments. In Section~\ref{subsec:theta_c} we study how the results are affected when we impose the arc length parametrisation and enforce the boundary conditions to be exactly satisfied by the network $\theta_{\boldsymbol{\rho}}^{\textrm{c}}$.

\subsection{Numerical experiments with \texorpdfstring{$q_{\boldsymbol{\rho}}^{\textrm{c}}$}{TEXT}}\label{subsec:num_exp_r_rho_c}

\begin{figure}[ht]
\centering
\includegraphics[width=0.48\textwidth]{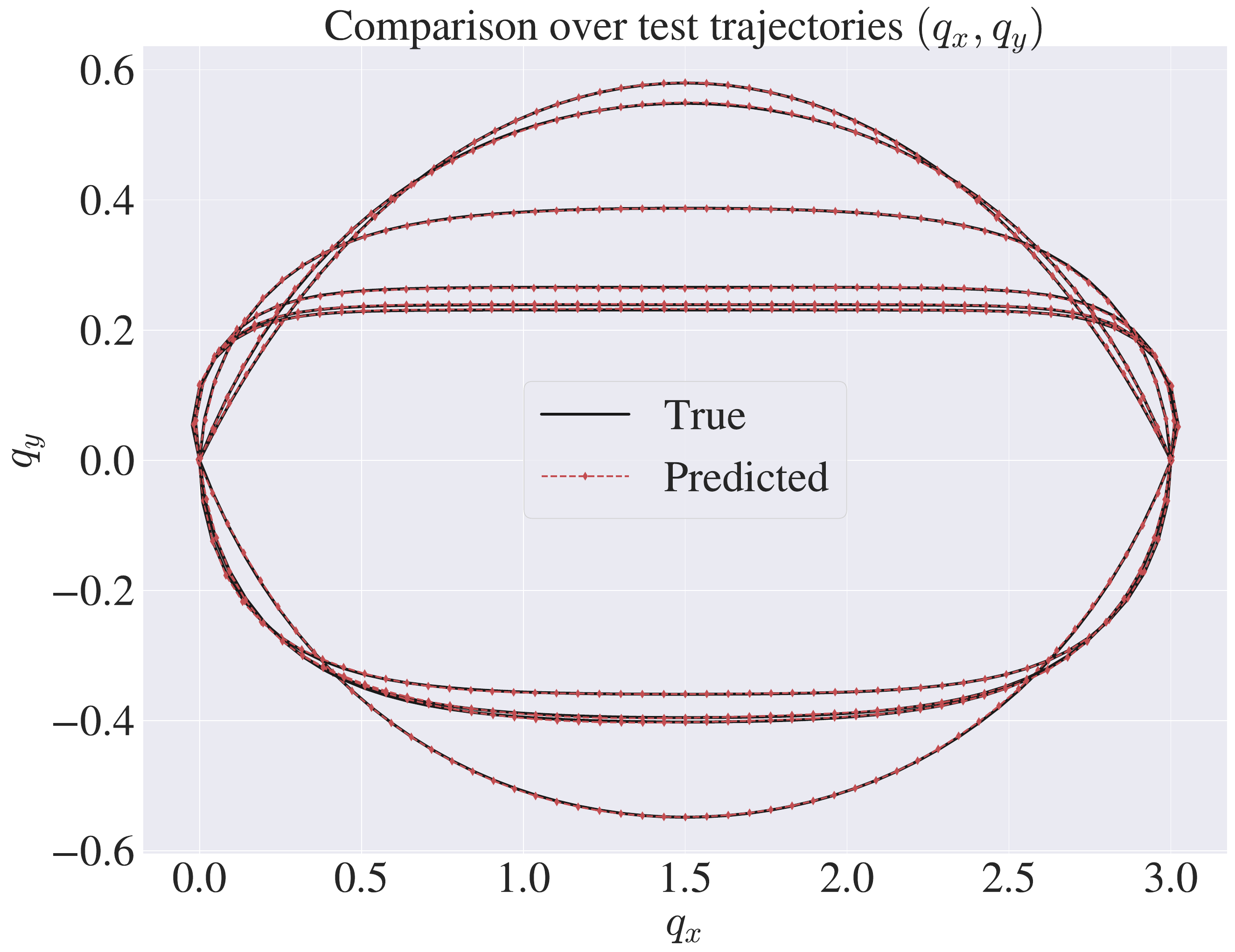}
\includegraphics[width=0.49\textwidth]{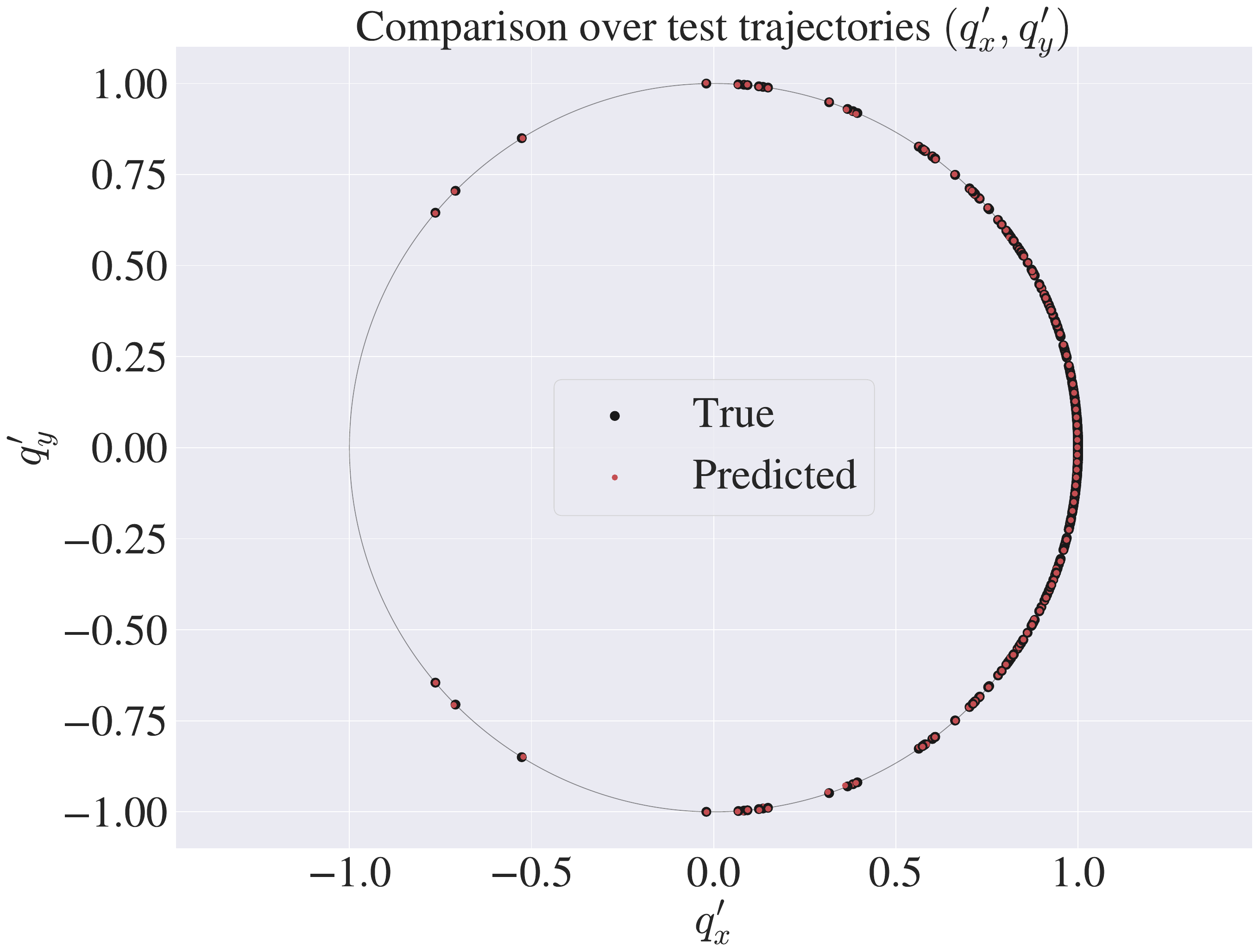}
\includegraphics[width = 0.49\textwidth]{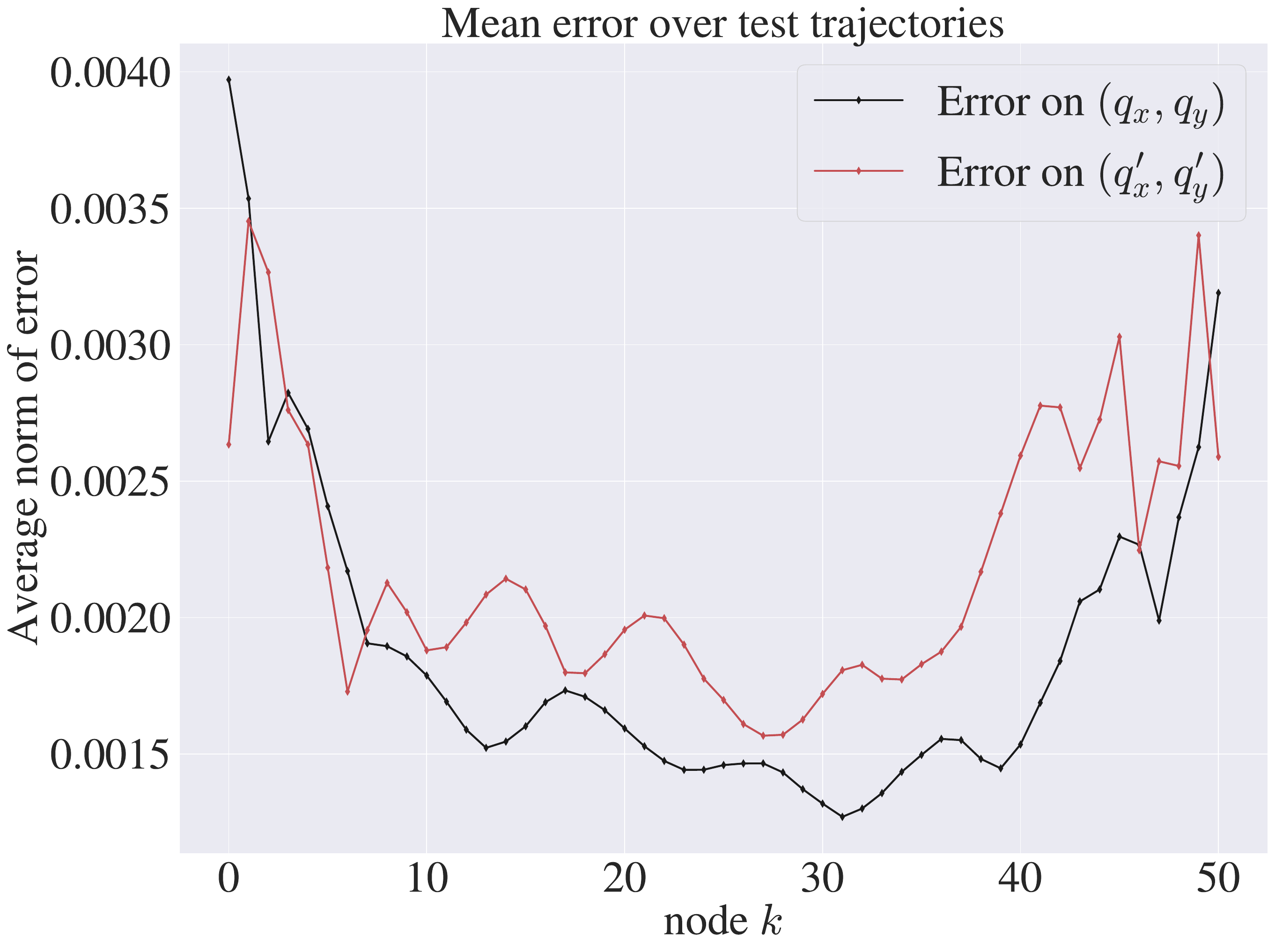}
\caption{Comparison over test trajectories for $\mathbf{q}$ and $\mathbf{q}^{\prime}$ for the continuous network $q_{\boldsymbol{\rho}}^{\textrm{c}}$ tested on the \textit{both-ends} data set with $80\% \-- 10\% \-- 10\%$ splitting into training, validation, and test sets. 
% \ergys{\sout{These results are obtained with the hyperparameters from Table \ref{tab:comp_table_hyperparameters_continuous_net}, that yield a training error equal to $3.554\cdot 10^{-6}$, a validation error equal to $4.779\cdot 10^{-6}$, and a test error equal to $4.354 \cdot 10^{-6}$.}} 
The mean squared error on the test set equals $4.354 \cdot 10^{-6}$. For presentation purposes, only 10 randomly selected trajectories are considered in the first two plots.}
\label{fig:q_qp_continuous_net}
\end{figure}

As for the case of the discrete network, we
perform an in-depth investigation of this learning setting. \newsub{In this case, the experiments are run on a GPU-P100 machine.} \newsub{For this continuous setup, the standard MLP architecture does not provide accurate results even after a hyperparameter optimisation routine. Given, hence, that the simple MLP architecture does not seem to be flexible enough to capture the complexity of the elastica solution in this continuous framework, we move to a different architecture that we call MULT for the presence of multiplicative interactions in its architecture. This network has demonstrated superior performance to standard fully connected neural networks in the context of operator learning, see e.g. \cite{wang2023long}. Details on this architecture can be found in Appendix \ref{sec:other_nns}.} \newsub{We fix the learning rate $\eta$ to $5\cdot 10^{-3}$ and only vary the number of layers and of hidden nodes in the training procedure, with the range of options reported in Appendix \ref{app2}, Table~\ref{tab:hyperparams_continuous_net}.} In this case, we define the loss as in Equation \eqref{eq:lossContinuos}, with $\gamma=10^{-2}$. The weight decay is systematically set to $0$. 
% \andout{Table~\ref{tab:best_hyperparams_cont_net} collects the combination of hyperparameters yielding the best results on the \textit{both-ends} data set.} 
\newsub{For the \textit{both-ends} data set, this leads to a training error equal to $  3.554 \cdot 10^{-6}$, a validation error equal to $ 4.779  \cdot 10^{-6}$, and a test error equal to $ 4.354 \cdot 10^{-6}$.} In Figure~\ref{fig:q_qp_continuous_net}, the comparison over test trajectories for $\mathbf{q}$ and $\mathbf{q}^{\prime}$ is shown.
As we can see in the plot showing the mean error over the trajectories, the error on the end nodes is nonzero, since we are not imposing boundary conditions by construction. This is in contrast to the corresponding plot for the discrete network in Figure~\ref{fig:q_qp_discrete_net}.

Also in this case, we examine the behaviour of the learning process with different splittings of the data set into training and test sets. We display the results in Table \ref{tab:comp_table_accuracy_continuous_net} and summarise the corresponding hyperparameters in Appendix \ref{app2}, Table \ref{tab:comp_table_hyperparameters_continuous_net}.
\begin{table*}[ht]
\small
\begin{tabularx}{1.\textwidth}{ 
   >{\centering\arraybackslash}p{0.28\textwidth} 
   >{\centering\arraybackslash}p{0.21\textwidth}   
   >{\centering\arraybackslash}p{0.21\textwidth}
   >{\centering\arraybackslash}p{0.21\textwidth}}
\toprule
$\begin{array}{c}
    \text{Data set splitting} \\
    \text{Training \-- validation \-- test}
    \end{array}$ & \text {Training accuracy} & \text {Validation accuracy} & \text{Test accuracy} \\
    \midrule \text {10\% \-- 10\% \-- 10\%} & $ 2.146  \cdot 10^{-4}$ & $ 1.252  \cdot 10^{-3}$ & $  8.811 \cdot 10^{-4}$\\
    \text {20\% \-- 10\% \-- 10\%} & $  4.187 \cdot 10^{-5}$  & $ 4.239  \cdot 10^{-5}$ & $ 6.279  \cdot 10^{-5}$\\
    \text {40\% \-- 10\% \-- 10\%} & $  7.037 \cdot 10^{-6}$ & $  8.357  \cdot 10^{-6}$ & $ 8.434 \cdot 10^{-6}$ \\
    \text {80\% \-- 10\% \-- 10\%} & $  3.554 \cdot 10^{-6}$ & $ 4.779  \cdot 10^{-6}$ & $ 4.354 \cdot 10^{-6}$\\
    \bottomrule
\end{tabularx}\caption{Behaviour of the continuous network $q_{\boldsymbol{\rho}}^{\textrm{c}}$ tested on the \textit{both-ends} data set with fewer training data points. The size of the training set varies, while that of the validation and the test sets is fixed. The last row corresponds to the results in Figure \ref{fig:q_qp_continuous_net}.}
\label{tab:comp_table_accuracy_continuous_net} 
\end{table*}
\subsection{Numerical experiments with \texorpdfstring{$\theta_{\boldsymbol{\rho}}^{\textrm{c}}$}{TEXT}}\label{subsec:theta_c}
\begin{figure}[ht]
\centering
\includegraphics[width=0.49\textwidth]{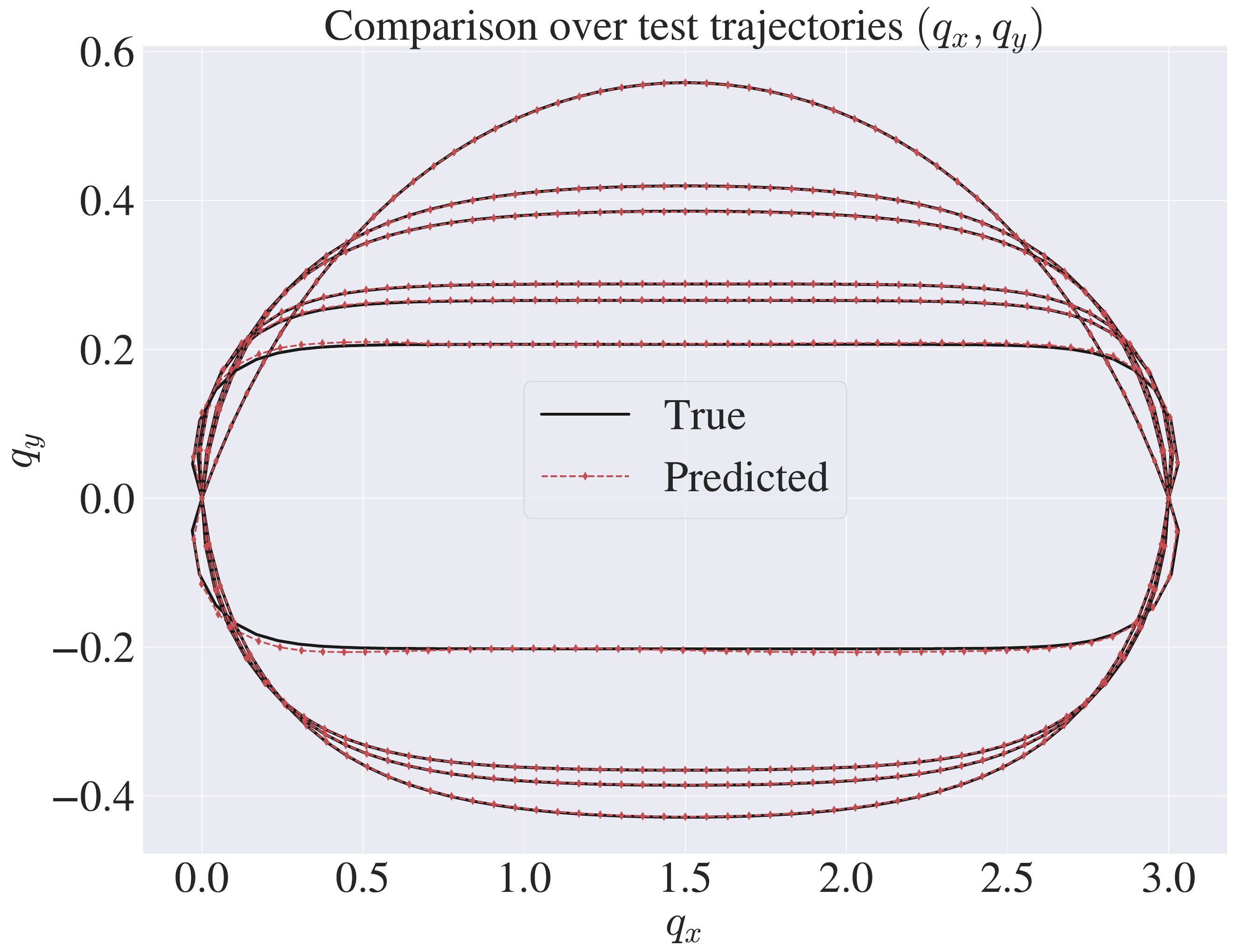} %q
\includegraphics[width=0.49\textwidth]{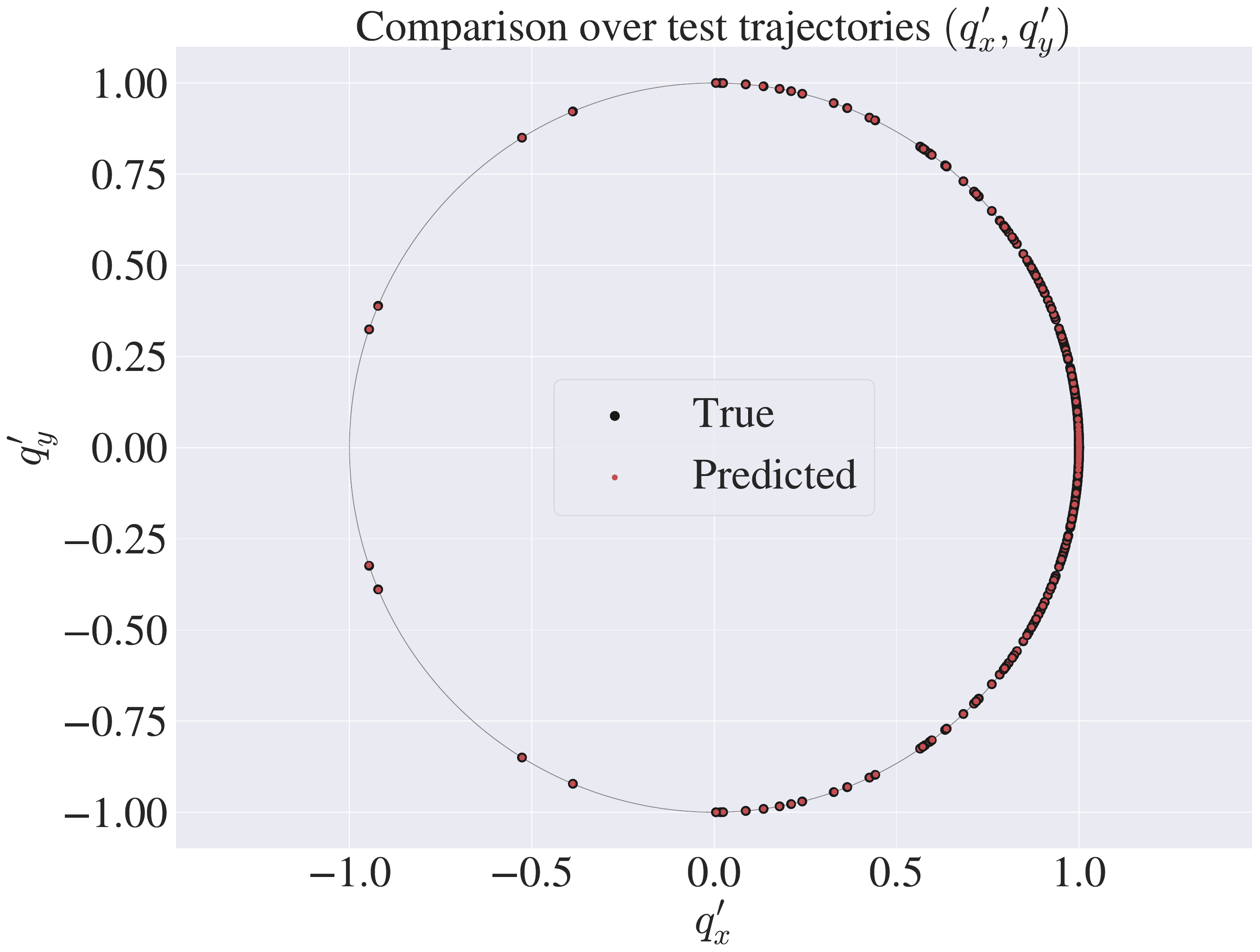} %v
\includegraphics[width=0.49\textwidth]{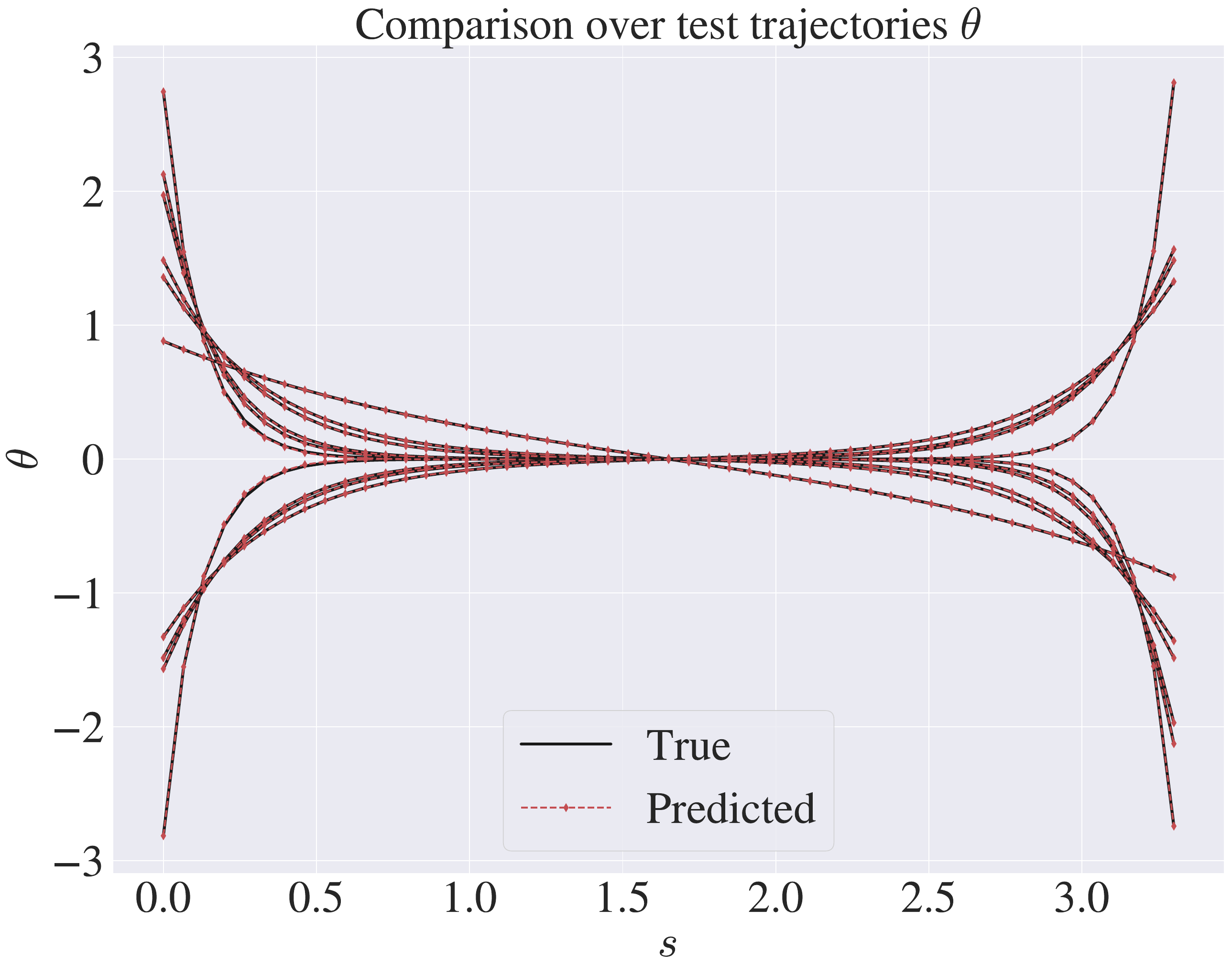} %theta
\includegraphics[width = 0.49\textwidth]{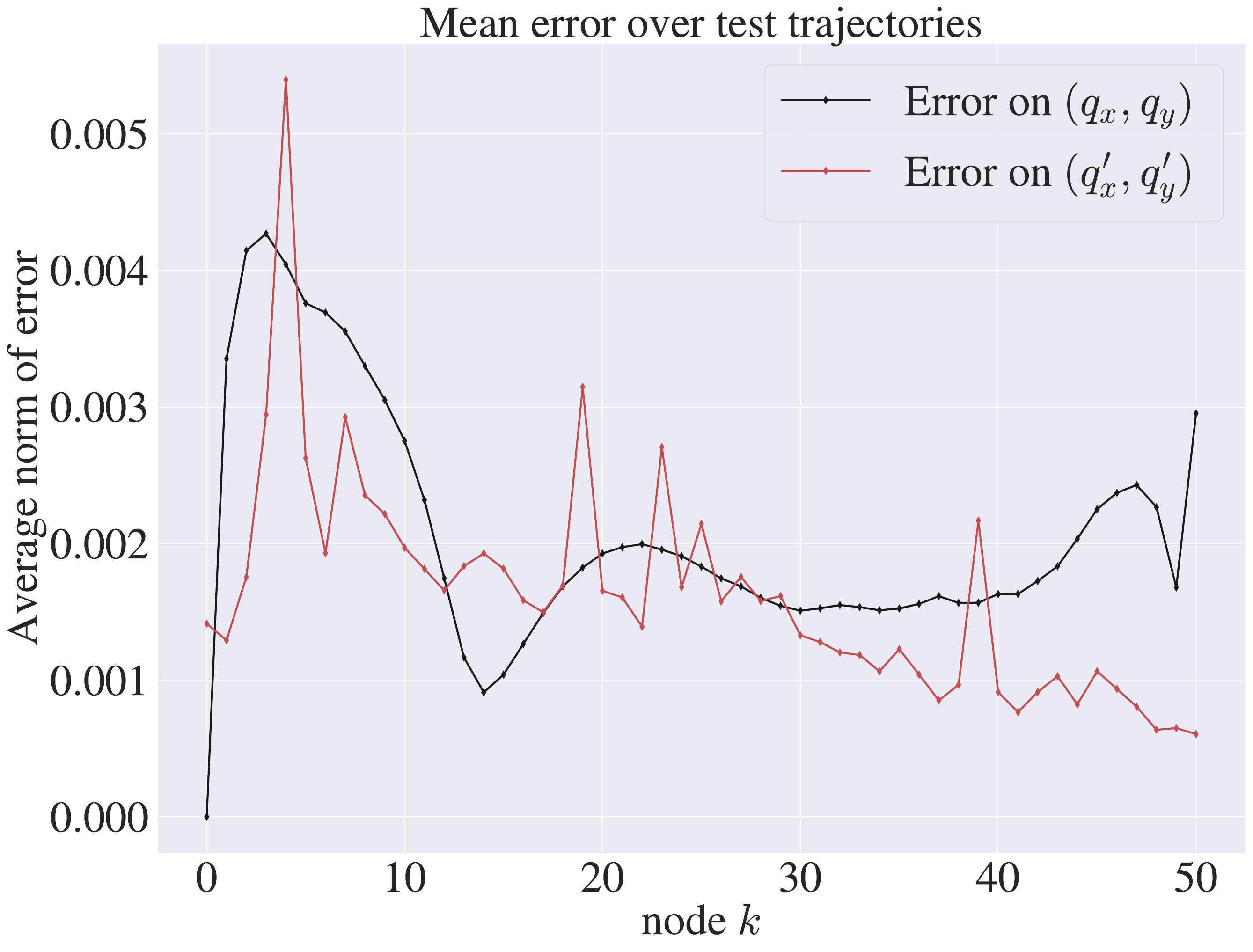} %error
\caption{Comparison over test trajectories for $\mathbf{q}$ and $\mathbf{q}^{\prime}$, for the case $\theta_{\boldsymbol{\rho}}^{\textrm{c}}$ is modelled as in Equation~\eqref{eq:noBCs}, with $80\% \-- 10\%\-- 10\%$ splitting of the \textit{both-ends} data set into training, validation, and test sets. 
% \ergys{\sout{These results are obtained with the hyperparameters from Table \ref{tab:hyperparams_selected_continuous_net_theta}, that yield a training error equal to $ 5.090\cdot 10^{-6}$, a validation error of $5.763\cdot 10^{-6}$, and a test error equal to $ 6.289 \cdot 10^{-6}$.}} 
The mean squared error on the test set equals $ 6.289 \cdot 10^{-6}$. For presentation purposes, only 10 randomly selected trajectories are considered in the first two plots.}
\label{fig:q_qp_continuous_net_theta_no_bc}
\end{figure}

\begin{figure}[ht]
\centering
\includegraphics[width=0.49\textwidth]{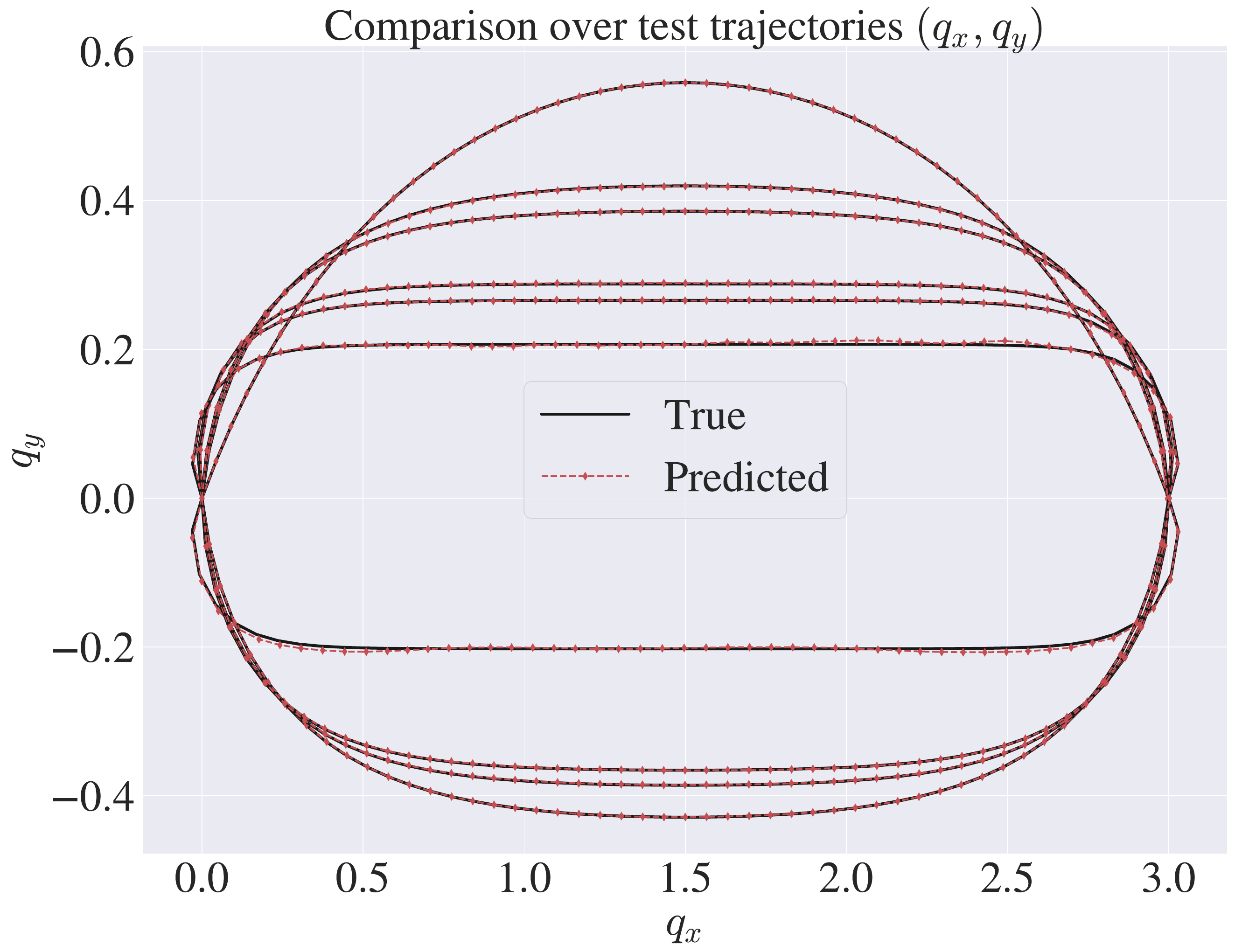} %q
\includegraphics[width=0.49\textwidth]{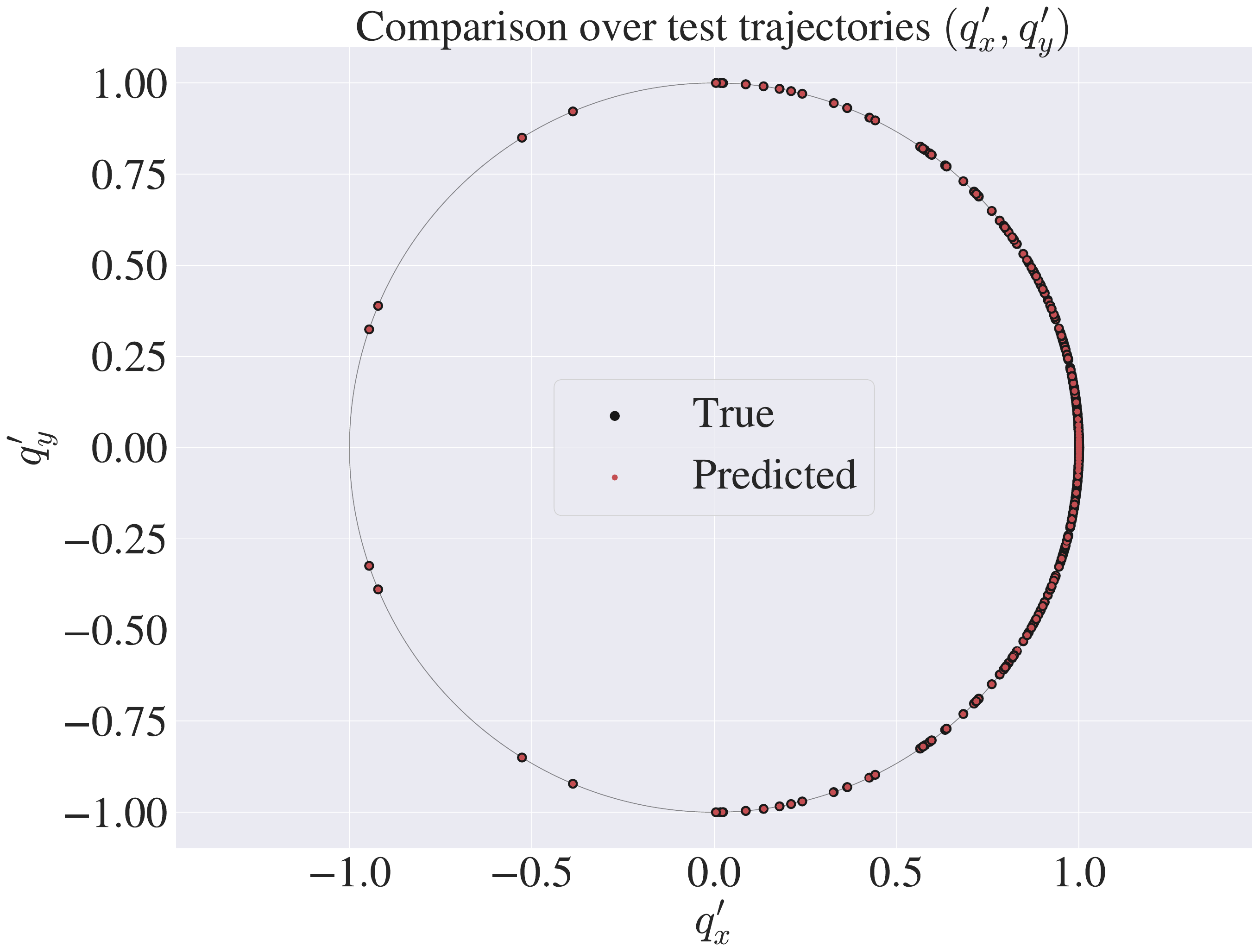} %v
\includegraphics[width = 0.49\textwidth]{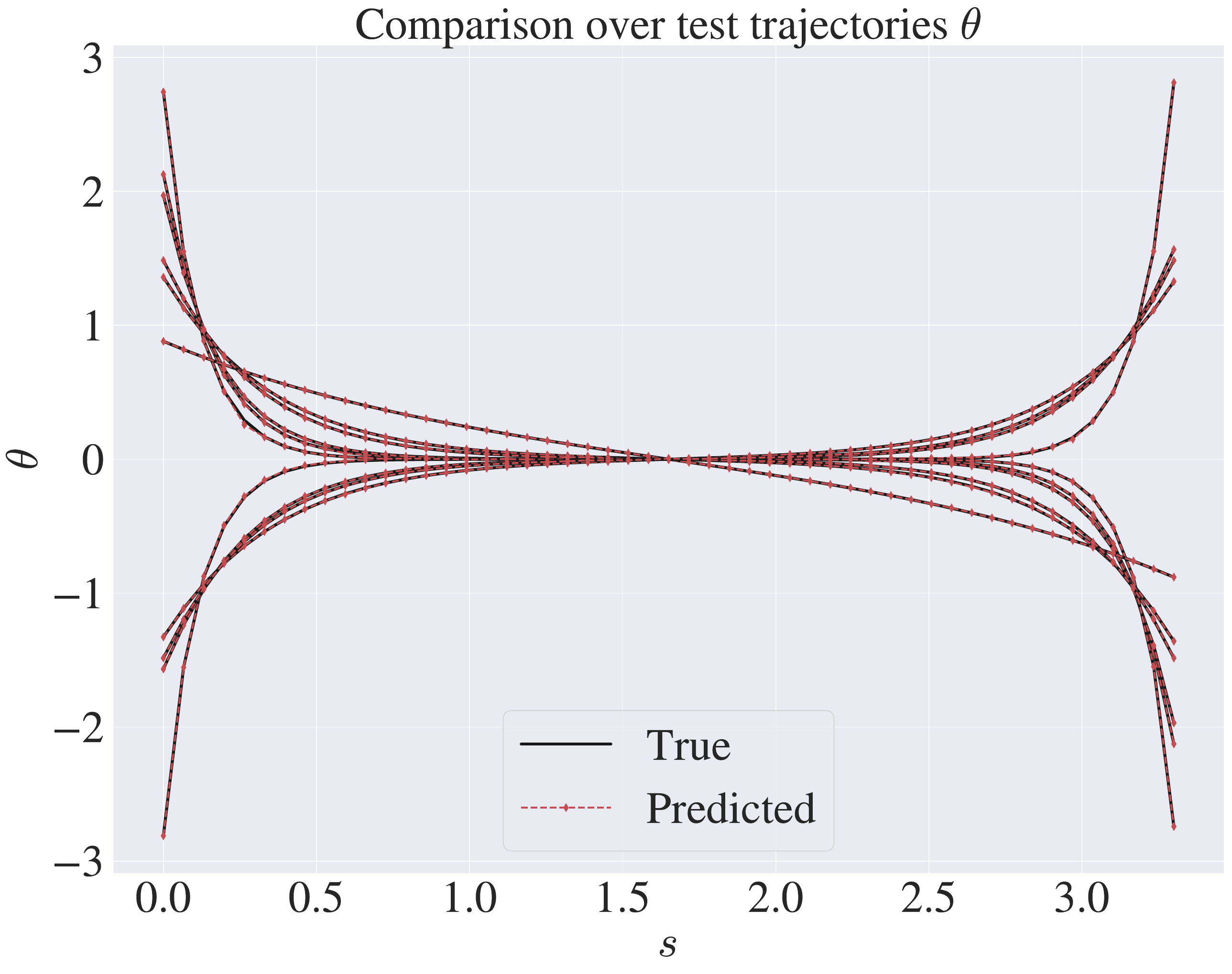} %theta
\includegraphics[width = 0.49\textwidth]{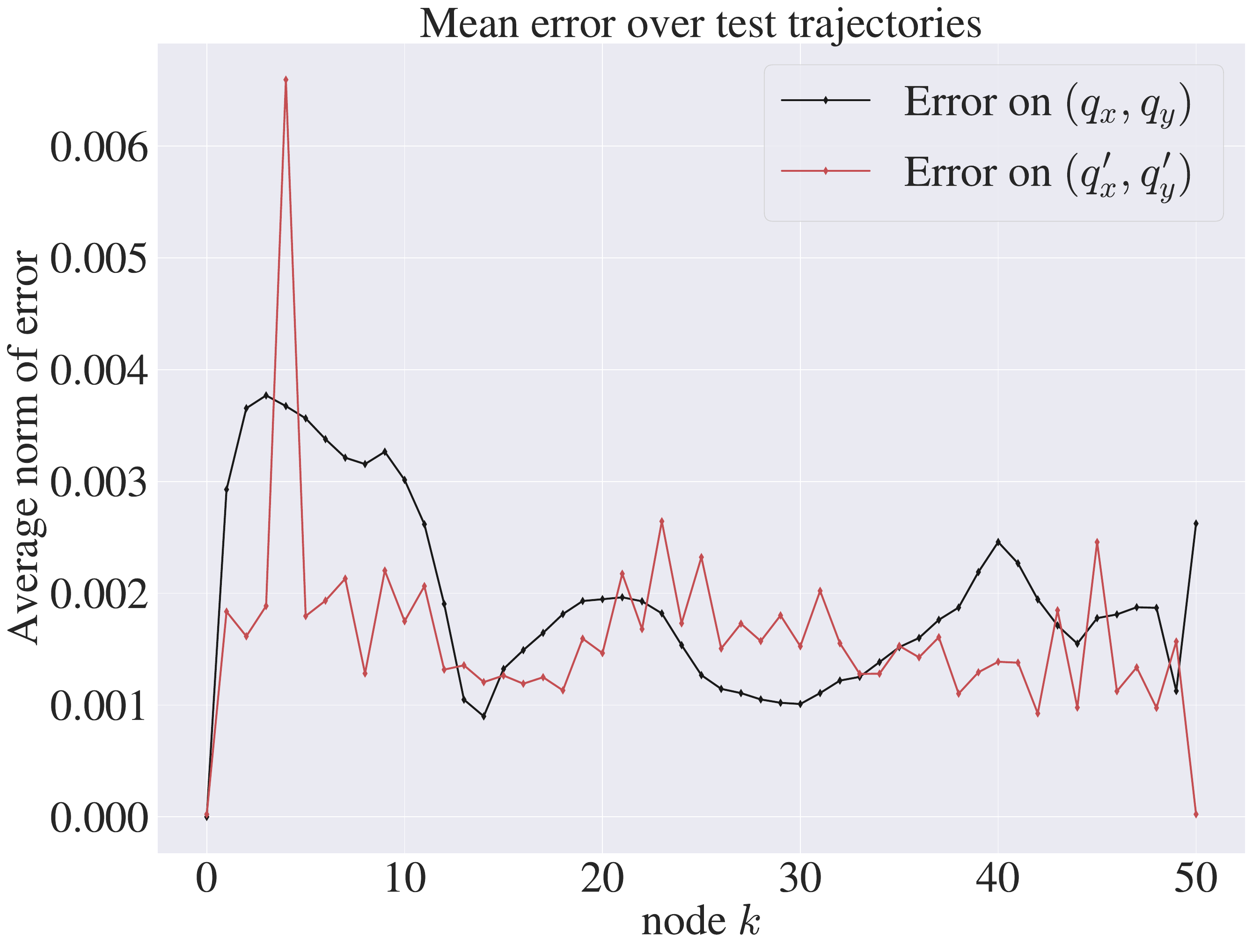} %error
\caption{Comparison over test trajectories for $\mathbf{q}$ and $\mathbf{q}^{\prime}$, for the case $\theta_{\boldsymbol{\rho}}^{\textrm{c}}$ is modelled as in Equation~\eqref{eq:withBCs}, with $80\% \-- 10\% \-- 10\%$ splitting of the \textit{both-ends} data set into training, validation, and test sets. 
% \ergys{\sout{These results are obtained with the hyperparameters from Table \ref{tab:hyperparams_selected_continuous_net_theta}, that yield a training error equal to $ 5.301\cdot 10^{-6}$, a validation error equal to $ 5.065\cdot 10^{-6}$, and a test error equal to $ 4.385 \cdot 10^{-6}$.}} 
The mean squared error on the test set equals $ 4.385 \cdot 10^{-6}$.
For presentation purposes, only 10 randomly selected trajectories are considered in the first two plots.}
\label{fig:q_qp_continuous_net_theta_bc}
\end{figure}

Here we consider a neural network approximation of the angle $\theta(s)$ that parametrises the tangent vector $\mathbf{q}'(s)=(\cos(\theta(s)),\sin(\theta(s)))$. By design, the approximation $\tau_{\boldsymbol{\rho}}^{\mathrm{c}}$ of the tangent vector $\mathbf{q}'$ satisfies the constraint $\|\tau_{\boldsymbol{\rho}}^{\mathrm{c}}(\mathbf{x}^i)(s)\|_2=1$ for every $s\in [0,L]$ and $\mathbf{x}^i\in\mathbb{R}^8$. We also analyse how the neural network approximation behaves when the boundary conditions $\tau_{\boldsymbol{\rho}}^{\mathrm{c}}(\mathbf{x}^i)(0)=\mathbf{q}'(0)$ and $\tau_{\boldsymbol{\rho}}^{\mathrm{c}}(\mathbf{x}^i)(L)=\mathbf{q}'(L)$ are imposed by construction. To do so, we model the parametric function $\hat{\theta}_{\boldsymbol{\rho}}^{\textrm{c}}$, defined in Equation \eqref{eq:thetaHat}, in one of the two following ways:
{
\begin{equation}\label{eq:noBCs}
\hat{\theta}_{\boldsymbol{\rho}}^{\textrm{c}}(\mathbf{x}^i)(s) = f_{\boldsymbol{\rho}}(s,\theta_0^i,\theta_N^i),
\end{equation}
}
{
\begin{equation}\label{eq:withBCs}
\begin{split}
    \hat{\theta}_{\boldsymbol{\rho}}^{\textrm{c}}(\mathbf{x}^i)(s) &= f_{\boldsymbol{\rho}}(s,\theta_0^i,\theta_N^i) + (\theta_0^i-f_{\boldsymbol{\rho}}(0,\theta_0^i,\theta_N^i))e^{-100s^2} \\
    &+ (\theta_N^i-f_{\boldsymbol{\rho}}(L,\theta_0^i,\theta_N^i))e^{-100(s-L)^2},
    \end{split}
\end{equation}
}
where $f_{\boldsymbol{\rho}}:\mathbb{R}^3\to\mathbb{R}$ is any neural network, and we recall that $\pi(\mathbf{x}^i)=(\theta_0^i,\theta_N^i)$. We remark that, in the case of the parameterisation in Equation~\eqref{eq:withBCs}, one gets $\theta_{\boldsymbol{\rho}}^{\textrm{c}}(\mathbf{x}^i)(0)=\theta_0^i$ and $\theta_{\boldsymbol{\rho}}^{\textrm{c}}(\mathbf{x}^i)(L)=\theta_N^i$ up to machine precision, due to the fast decay of the Gaussian function. As in the previous sections, we collect the hyperparameter and architecture options with the respective range of choices in Table \ref{tab:hyperparams_continuous_net_theta}, and we report the results without imposing the boundary conditions in Figure \ref{fig:q_qp_continuous_net_theta_no_bc}, while those imposing them in Figure \ref{fig:q_qp_continuous_net_theta_bc}, in both cases using the \textit{both-ends} data set, with $80\% \-- 10\% \-- 10\%$ splitting into training, validation, and test sets. The results shown in the two figures correspond respectively to training errors of $ 6.288 \cdot 10^{-6}$ and $ 5.301\cdot 10^{-6}$, validation errors $ 5.874 \cdot 10^{-6}$ and $ 5.065 \cdot 10^{-6}$, and test errors of $ 5.089 \cdot 10^{-6}$ and $ 4.385 \cdot 10^{-6}$. The best-performing hyperparameter combinations can be found in Table \ref{tab:hyperparams_selected_continuous_net_theta}.

\section{Discussion}\label{conclusion}

The results in Figures \ref{fig:q_qp_continuous_net_theta_no_bc} and \ref{fig:q_qp_continuous_net_theta_bc} are comparable, especially looking at the mean error plots. This suggests that the imposition of the boundary conditions, in the proposed way, is not limiting the expressivity of the considered network. Thus, given the boundary value nature of our problem, these figures advocate the enforcement of the boundary conditions on the network $\theta_{\boldsymbol{\rho}}^{\textrm{c}}$. However, due to the chosen reconstruction procedure in Equation~\eqref{eq:q_rec} for the variable $\mathbf{q}$, we are able to impose the boundary conditions on $\mathbf{q}$ only on the left node. Other more symmetric reconstruction procedures can be adopted, but the proposed one has provided better experimental results.

Comparing the results related to $q_{\boldsymbol{\rho}}^{\textrm{c}}$ with those of $\theta_{\boldsymbol{\rho}}^{\textrm{c}}$, we notice similar performances in terms of training and test errors. In both cases, they have one order of magnitude more than the corresponding training and test errors of the discrete network $q_{\boldsymbol{\rho}}^{\textrm{d}}$. Thus, as a result of our experiments, we can conclude that
\begin{itemize}
    \item if the accuracy and the efficient evaluation of the model at the discrete nodes are of interest, the discrete network is the best option;
    \item for a more flexible model, not restricted to the discrete nodes, the continuous network is a better choice; among the two proposed modelling strategies, working with $q_{\boldsymbol{\rho}}^{\textrm{c}}$ is more suitable for an easy parametrisation of both $\mathbf{q}$ and $\mathbf{q}'$, while $\theta_{\boldsymbol{\rho}}^{\textrm{c}}$ is more suitable to impose geometrical structure and constraints. 
\end{itemize}
The total accuracy error of a neural network model can be defined by splitting it into three components: approximation error, optimisation error, and generalisation error (see e.g. \cite{lu2021learning}). To achieve excellent agreement between predicted and reference trajectories, it is crucial to select the appropriate architecture and fine-tune the model hyperparameters. Our results demonstrate that we can construct a network that is expressive enough to provide a small approximation error and with very good generalisation capability. 

\newsub{Lastly, we have compared the time cost of the Neural Network prediction against the traditional approach with numerical solvers as described in Section \ref{subsec:bifurcation}. The discrete and the continuous approaches outperformed the traditional solvers with an average speedup of 105.000 times and 260.000 times, respectively, across the test trajectories. The training time of the continuous network is 1.25 times larger than that of the discrete network. It's important to note that these results are subject to certain limitations, such as the specific choice of the hyperparameters or the machine used to train and test the network.
These findings suggest that using neural networks to predict new solutions of the elastica for unseen boundary conditions is much more time efficient than the classical numerical methods, although requiring intensive offline training.}

\subsection{Future work} \label{pinn}

In the methods presented in this paper, the mathematical problem and the neural network model do not interact once the data set is created. To improve the results presented here, one could include Euler's elastica model directly into the training process.
This could be done either by directly imposing in the loss function that $\mathbf{q}(s)$ satisfies the differential equations \eqref{eq:second_order_ELeq}, or one could add the constrained action integral from Equation \eqref{eq:action_integral} into the loss function that is minimised, see e.g. \cite{lagaris1998artificial, raissi2019physics, yu2018deep, samaniego2020energy}.

There are many promising directions to follow up on this work. One is to consider 3D versions of Euler's elastica, another is to look at the dynamical problem, and finally one may examine industrial applications, as mentioned in the introduction.

% \bigskip

% \textbf{Acknowledgments.} This project has received funding from the European Union’s Horizon 2020 research and innovation programme under the Marie Skłodowska-Curie grant agreement No 860124. This work was partially supported by a grant from the Simons Foundation (DM). This contribution reflects only the authors’ view, and the Research Executive Agency and the European Commission are not responsible for any use that may be made of the information it contains. 

\bigskip

\section*{CRediT author statement.} \textbf{Elena Celledoni:} Conceptualisation, Validation, Writing - Review \& Editing, Supervision, Funding acquisition. \textbf{Ergys \c{C}okaj:} Validation, Investigation, Visualisation, Writing - Review \& Editing. \textbf{Andrea Leone:} Methodology, Software, Investigation, Writing - Original Draft, Review \& Editing. \textbf{Sigrid Leyendecker:} Conceptualisation, Writing - Review \& Editing, Supervision, Funding acquisition. \textbf{Davide Murari:} Methodology, Software, Investigation, Writing - Original Draft, Review \& Editing. \textbf{Brynjulf Owren:} Conceptualisation, Validation, Writing - Review \& Editing, Supervision, Funding acquisition. \textbf{Rodrigo T. Sato Mart\'{i}n de Almagro:} 
Methodology, Validation, Writing - Review \& Editing, Supervision. \textbf{Martina Stavole:} Methodology, Software, Investigation, Writing - Original Draft, Review \& Editing.

\bigskip

\textbf{Acknowledgments.} This project has received funding from the European Union’s Horizon 2020 research and innovation programme under the Marie Skłodowska-Curie grant agreement No 860124. This work was partially supported by a grant from the Simons Foundation (DM). This contribution reflects only the authors’ view, and the Research Executive Agency and the European Commission are not responsible for any use that may be made of the information it contains. 

\appendix

\section{\newsub{Architecture for the continuous network}}\label{sec:other_nns}
We provide the expression of the forward propagation of the multiplicative network MULT used for the experiments in Section~\ref{sec:nncont}:
{
\begin{align} %\label{eq:multiplicativenn}
&\mathbf{U}=\sigma(\mathbf{W}_1\mathbf{x}+\mathbf{b}_1), \quad \mathbf{V}=\sigma(\mathbf{W}_2\mathbf{x} +\mathbf{b}_2)\label{eq:fistStepMult}\\
&\mathbf{H}_1 = \sigma(\mathbf{W}
_3\mathbf{x} +\mathbf{b}_3)\\
&\mathbf{Z}_{j} = \sigma(\mathbf{W}^{z}_j\mathbf{H}_j + \mathbf{b}^{z}_j),\;j=1,\dots,\ell\\
&\mathbf{H}_{j+1}=(1-\mathbf{Z}_j)\odot \mathbf{U}+ \mathbf{Z}_j\odot \mathbf{V},\;j=1,\dots,\ell\\
&f_{\boldsymbol{\rho}}^{\mathrm{MULT}}(\mathbf{x}) = \mathbf{W}\mathbf{H}_{\ell+1}+\mathbf{b}\label{eq:lastStepMult},
\end{align}
}
where $\odot$ denotes the component-wise multiplications. In this case, $\boldsymbol{\rho} = \left\{ \mathbf{W}_1, \mathbf{b}_1, \mathbf{W}_2, \mathbf{b}_2, \mathbf{W}_3, \mathbf{b}_3, \left(\mathbf{W}^z_j, \mathbf{b}^z_j\right)_{j=1}^{\ell}, \mathbf{W}, \mathbf{b} \right\}$, and the weight matrices and biases have shapes that allow for the expressions \eqref{eq:fistStepMult}-\eqref{eq:lastStepMult} to be well-defined. This architecture is inspired by neural attention mechanisms and was introduced in \cite{wang2021understanding} to improve the gradient behaviour. A further motivation for our choice of including this architecture is experimental since it has proven effective in solving the task of interest, while still having a similar number of parameters to the MLP architecture. Throughout the paper, we refer to this architecture as \textit{multiplicative} since it includes component-wise multiplications, which help capture multiplicative interactions between the variables.

\section{\newsub{Details on hyperparameter optimisation}} \label{app2}

\newsub{We provide here further details on the hyperparameter optimisation strategy with \texttt{Optuna}, relative to the results in Sections \ref{nndiscrete} and \ref{sec:nncont}.} \newsub{The tables below display the hyperparameters we optimise for in each of the networks, the ranges and the distribution, as well as the selected combinations used to perform the experiments in the paper.}
\begin{table*}[ht]
\centering
\small
{
\begin{tabularx}{1.\textwidth}{ 
   >{\centering\arraybackslash}p{0.30\textwidth} 
   >{\centering\arraybackslash}p{0.30\textwidth} 
   >{\centering\arraybackslash}p{0.30\textwidth}}

\toprule
Hyperparameter  & Range & Distribution\\
\midrule
     \#layers $\ell$ & $\{0, ...,10\}$ & discrete uniform\\
     \#hidden nodes &  $[10,1000]\cap \mathbb{N}$ & discrete uniform\\
     $\gamma$ &  $[0, 1 \cdot 10^{-2}]$ & uniform\\
 \bottomrule
\end{tabularx}
}
    \caption{\newsub{Hyperparameter ranges for the discrete network $q_{\boldsymbol{\rho}}^{\textrm{d}}$. The first column of the table reports the hyperparameters we test. The second describes the set of allowed values for each, while the third specifies how such values are explored through \texttt{Optuna}.}}
\label{tab:hyperparams_dicrete_net}
\end{table*}
\begin{table*}[ht]
\small
\begin{tabularx}{1.\textwidth}{ 
   >{\centering\arraybackslash}p{0.3\textwidth}
   >{\centering\arraybackslash}p{0.135\textwidth} 
   >{\centering\arraybackslash}p{0.135\textwidth} 
   >{\centering\arraybackslash}p{0.135\textwidth} 
   >{\centering\arraybackslash}p{0.135\textwidth}}
\toprule
\multirow{2}{*}{$\begin{array}{c}
    \text{Hyperparameter} \\
    \text {combination}
    \end{array}$}&\multicolumn{4}{c}{\newsub{\% of trajectories of the whole dataset in the training set}}\\ 
    \cmidrule{2-5}
    & {10\%} & {20\%} & {40\%} & {80\%}\\
\midrule
 \text{\# layers $\ell$} & 4 & 4 & 4 & 4 \\
    \text{\#hidden nodes}  & 950 & 978 & 997 & 985 \\
 $\gamma$  & $ 7.044\cdot 10^{-3}$ & $ 6.336\cdot 10^{-3}$ &  $ 9.004\cdot 10^{-3}$ &  $3.853\cdot 10^{-3}$ \\
\bottomrule
\end{tabularx}\caption{\label{tab:comp_table_hyperparameters_discrete_net}\newsub{Choice of hyperparameters for the training of the discrete network $q_{\boldsymbol{\rho}}^{\textrm{d}}$ tested on the \textit{both-ends} data set with different sizes of the training set, with the validation and test sets each containing $10\%$ of trajectories of the dataset.}}
\end{table*}
\begin{table*}[ht]
\centering
\small
{
\begin{tabularx}{1.\textwidth}{ 
   >{\centering\arraybackslash}p{0.30\textwidth}
   >{\centering\arraybackslash}p{0.30\textwidth}
   >{\centering\arraybackslash}p{0.30\textwidth}}
\toprule
Hyperparameter   & Range & Distribution\\
\midrule
     \#layers $\ell$ & $\{5, \ldots, 10\}$ & discrete uniform\\
     \#hidden nodes &  $[10, 250] \cap \mathbb{N} $ & discrete uniform\\
 \bottomrule
\end{tabularx}
}
\caption{\newsub{Hyperparameter ranges for the continuous network $q_{\boldsymbol{\rho}}^{\textrm{c}}$. The first column of the table reports the hyperparameters we test. The second describes the set of allowed values for each, while the third specifies how such values are explored through \texttt{Optuna}.}
% Hyperparameter ranges for the continuous network $q_{\boldsymbol{\rho}}^{\textrm{c}}$ tested on the \textit{both-ends} data set with $80\% \-- 10\% \-- 10\%$ splitting into training, validation, and test sets. The first column of the table reports the hyperparameters and network architectures we test for. The second describes the set of allowed values for each, while the third specifies how such values are explored through \texttt{Optuna}. Finally, the fourth column shows the combination of hyperparameters yielding the best result, corresponding to Figure \ref{fig:q_qp_continuous_net}. \andout{The weight decay is systematically set to $0$. MULT stands for multiplicative neural network and corresponds to the network in Equations~\eqref{eq:fistStepMult}-\eqref{eq:lastStepMult} of Appendix \ref{sec:other_nns}.}
}
\label{tab:hyperparams_continuous_net}
\end{table*}
\begin{table*}[ht]
\small
\begin{tabularx}{1.\textwidth}{ 
   >{\centering\arraybackslash}p{0.3\textwidth}
   >{\centering\arraybackslash}p{0.135\textwidth} 
   >{\centering\arraybackslash}p{0.135\textwidth} 
   >{\centering\arraybackslash}p{0.135\textwidth} 
   >{\centering\arraybackslash}p{0.135\textwidth}}
\toprule
\multirow{2}{*}{$\begin{array}{c}
    \text{Hyperparameter} \\
    \text {combination}
    \end{array}$}&\multicolumn{4}{c}{\newsub{\% of trajectories of the whole dataset in the training set}}\\ 
    \cmidrule{2-5}
    & {10\%} & {20\%} & {40\%} & {80\% }\\
\midrule
 \text{\# layers $\ell$} & 6 & 7 & 8 &  6 \\
 \text{\#hidden nodes}  & 139  & 185 & 181 & 106  \\
\bottomrule
\end{tabularx}\caption{\newsub{Choice of hyperparameters for the training of the continuous network $q_{\boldsymbol{\rho}}^{\textrm{c}}$ tested on the \textit{both-ends} data set with different sizes of the training set, with the validation and test sets each containing $10\%$ of trajectories of the dataset.}}
\label{tab:comp_table_hyperparameters_continuous_net}
\end{table*}
\begin{table*}[ht]
\centering
\small
{
\begin{tabularx}{1.\textwidth}{ 
   >{\centering\arraybackslash}p{0.30\textwidth}
   >{\centering\arraybackslash}p{0.30\textwidth}
   >{\centering\arraybackslash}p{0.30\textwidth}}
\toprule
Hyperparameter   & Range & Distribution\\
\midrule
     \#layers $\ell$ & $\{1, \ldots, 10\}$ & discrete uniform\\
     \#hidden nodes &  $[50, 200] \cap \mathbb{N} $ & discrete uniform\\
 \bottomrule
\end{tabularx}
}
\caption{\newsub{Hyperparameter ranges for the continuous network $\theta_{\boldsymbol{\rho}}^{\textrm{c}}$.The first column of the table
reports the hyperparameters we test. The second describes the set of allowed values for each, while the third specifies how such values are explored through \texttt{Optuna}.}}
\label{tab:hyperparams_continuous_net_theta}
\end{table*}
\begin{table*}[ht]
\centering
\small
\begin{tabularx}{1.\textwidth}{ 
   >{\centering\arraybackslash}p{0.3\textwidth}
   >{\centering\arraybackslash}p{0.3\textwidth}
   >{\centering\arraybackslash}p{0.3\textwidth}}
\toprule
\multirow{2}{*}{$\begin{array}{c}
    \text{Hyperparameter} \\
    \text {combination}
    \end{array}$}&\multicolumn{2}{c}{{ $\theta_{\boldsymbol{\rho}}^{\textrm{c}}$}}\\ 
    \cmidrule{2-3}
    & $\theta_{\boldsymbol{\rho}}^{\textrm{c}}$ as in \eqref{eq:noBCs}  & $\theta_{\boldsymbol{\rho}}^{\textrm{c}}$ as in \eqref{eq:withBCs}\\
\midrule
     \# layers $\ell$ & 8 & 8\\
     \#hidden nodes & 93 & 58\\
 \bottomrule
\end{tabularx}
\caption{{Choice of the hyperparameters for the training of the continuous network $\theta_{\boldsymbol{\rho}}^{\textrm{c}}$ tested on the both-ends data set with $80\% \-- 10\%\-- 10\%$ splitting into training, validation, and test sets. The second column shows the combination of hyperparameters yielding the best result corresponding to Figure \ref{fig:q_qp_continuous_net_theta_no_bc}, while the third column that corresponding to Figure \ref{fig:q_qp_continuous_net_theta_bc}.}}
\label{tab:hyperparams_selected_continuous_net_theta}
\end{table*}

\bibliography{referencesV2} 

\end{document}